\documentclass[journal, onecolumn, draftcls, 12pt]{IEEEtran}
\usepackage{hyperref}
\usepackage{graphics, subfigure,  times, amsfonts, amsmath, amssymb, cite}
\usepackage{tikz}
\usepackage{url}
\usetikzlibrary{shapes,arrows}
\usepackage{pgfplots}
\usepackage{color}
\usepackage{multirow}
\usetikzlibrary{automata,positioning}
\usepackage{rotating}
\usepackage{flushend}
\usepackage{amsthm}
\usepackage{epstopdf}
\newtheorem{theorem}{Theorem}
\newtheorem{assumption}{Assumption}
\newtheorem{lemma}{Lemma}
\newtheorem{definition}{Definition}

\newtheorem{corollary}{Corollary}
\usepackage{bbm}
\title{ Control of Charging of Electric Vehicles through Menu-Based Pricing}
\date{}
\author{Arnob Ghosh and Vaneet Aggarwal,~\IEEEmembership{Senior~Member,~IEEE}\thanks{The authors are with the School of  Industrial Engineering, Purdue University, West Lafayette, IN 47907, email: \{ghosh39,vaneet\}@purdue.edu. This work was supported in part by the U.S. National Science Foundation
		under grant CCF-1527486.
		
		This paper was presented in part at the IEEE International Conference on Communications, Paris, France, May 2017 \cite{icc_veh17}. 
		
}}
\begin{document}
\maketitle
\begin{abstract}
We propose an online pricing mechanism for electric vehicle (EV) charging. A charging station decides prices for each arriving EV depending on the energy and the time within which the EV will be served (i.e. deadline).  The user selects either one of the contracts by paying the prescribed price or rejects all depending on their surpluses. The charging station can serve users using renewable energy and conventional energy. Users  may select longer deadlines as they may have to pay less because of the less amount of conventional energy, however, they have to wait a longer period.  We consider a {\em myopic} charging station and show that there exists a pricing mechanism which jointly maximizes the social welfare and the profit of the charging station when the charging station knows the utilities of the users. However, when the charging station does not know the utilities of the users, the social welfare pricing strategy may not maximize the expected profit of the charging station and even the profit may be $0$.  We propose a fixed profit pricing strategy which provides a guaranteed fixed profit to the charging station and can maximize the profit in practice.  We  empirically show that our proposed mechanism reduces the peak-demand and utilizes the limited charging spots in a charging station efficiently.
\end{abstract}


\section{Introduction}
Electric Vehicles (EVs) have several advantages over the traditional gasoline powered vehicles. For example, EVs are more environment friendly and more energy efficient. Thus, regulators (e.g. Federal Energy Regulator Commission (FERC)) are providing incentives to the consumers to switch to electric vehicles. However, the successful deployment of charging stations crucially depends on the profit of the charging stations and how efficiently the  resources are used for charging the electric vehicles.   Without profitable charging stations, the wide deployment of the EVs will remain a distant dream. On the other hand, because of the environment-friendly nature of the electric vehicles, it is also important for the regulators to  increase the user (or, consumer) surplus to provide an incentive for the users. Hence,  selecting a price is an imperative issue for the charging stations. The charging station may have limited charging spots or renewable energy harvesting devices. Hence, intelligent allocation of the resources among the EVs is a key component for fulfilling the potential of EVs' deployment.

We propose a menu based pricing scheme for charging an EV. Whenever an EV arrives, the charging station offers  a variety of  contracts $(l,t_{dead})$ at price $p_{l,t_{dead}}$  to the user where the user will be able to use  up to $l$ units of energy within the deadline $t_{dead}$  for completion. The EV user either accepts one of the contracts by paying the specified price or rejects all of those based on its payoff. We assume that the user gets an utility for consuming $l$ amount of energy with the deadline $t_{dead}$. The payoff of the user (or, user's surplus) for a contract is the difference between the utility and the price paid for one contract. The user will select the option which fetches the highest payoff. 

The various advantages of the above pricing scheme should be noted. First, it is an online pricing scheme. It can be adapted for each arriving user.  Second, since the charging station offers prices for different levels of energy and the deadline, the charging station can prioritize one contract over the others depending on the energy resources available. Favorable prices for shorter deadlines can attract users to vacate the charging stations early and only use it when it is necessary.  Third, the user's decision is much simplified.  She only needs to select one of the contracts (or, reject all) and will receive the prescribed amount within the prescribed deadline. Fourth, the pricing mechanism is inherently {\em individual rational}\footnote{Individual rationality means that the user will only select one of the contracts if it gets non-negative payoff.}, {\em incentive compatible}\footnote{Individual compatibility denotes that the user will select only that contract which gives the optimal payoff.}, and {\em truthful}\footnote{Truthfulness implies that the users will not achieve a higher payoff  by delaying their arrivals or selecting a sub-optimal strategy. Hence, the pricing mechanism is robust against the strategy selection of the users.}.

We consider that the charging station is equipped with renewable energy harvesting devices and a storage device for storing energies. The charging station may also  buy conventional energy from the market to fulfill the contract of the user if required. Hence, if a new user accepts  the contract $(l,t_{dead})$, a cost is incurred to the charging station. This cost may also depend on the existing EVs and their resource requirements. Hence, the charging station needs to find the optimal cost for fulfilling each contract. We show that obtaining that cost is equivalent to solve a {\em linear programming} problem.


We consider two optimization problems--i) social welfare\footnote{Social welfare is the sum of the profit of the charging station and the user surplus.} maximization, and  ii) the EV charging station's profit maximization. We first propose a pricing scheme  which maximizes the social welfare {\em irrespective of whether the charging station is aware of the utilities of the users or not}. The pricing scheme is  simple to compute, as the charging station selects a price which is equal to the marginal cost for fulfilling a certain contract for a new user (Theorem~\ref{thm:pricestrategy}, Corollary \ref{thm:price_expected}).   However, the above pricing scheme only provides {\em zero} profit to the charging stations. Thus, such a pricing scheme may not be useful to the charging station. We show that when a charging station is {\em clairvoyant} (i.e., the charging station knows the utilities of the users), there exists a pricing scheme which satisfies both the objectives (Theorem~\ref{thm:profit_max}). Though  in the above pricing mechanism, the user's surplus becomes $0$. Thus, a {\em clairvoyant} charging station may not be beneficial for the user's surplus.

The charging station may not know the exact utilities of the users, however, it may know the distribution function\footnote{We do not put any assumption on the distribution function} from where it is drawn. We investigate the existence of a pricing mechanism which maximizes the ex-post social welfare i.e. maximizes the social welfare  for every possible realization of the utility function.  In the scenario where the charging station does not know the exact utilities of the users, we show that there {\em may not exist a pricing strategy which simultaneously maximizes the ex-post social welfare and  the expected profit}. One has to give away the ex-post social welfare maximization in order to achieve expected profit maximization.  Thus, unlike when the charging station {\em is clairvoyant}, there may not exist a pricing strategy which simultaneously satisfies both the objectives when the exact utilities are unknown.  We propose a pricing strategy which can fetch the highest possible profit to the charging station under the condition that it  maximizes the ex-post social welfare (Theorem~\ref{thm:profitmax_uncertainty}).   Above pricing strategy provides a {\em worst case} maximum profit to the charging station.  We show that such a pricing strategy can fetch a higher profit when the charging station can harvest a large amount of renewable energy. However, the profit only increases up to a certain threshold, beyond that threshold harvested energy has no effect on the profit. 

Since the above pricing strategy  may not yield the {\em maximum expected} profit to the charging station, we have to relax the constraint the social welfare to be maximized in order to yield a higher profit to the charging station. Whether a contract will be selected by the user does not depend on the price of the contract, but also the prices of other contracts. Thus, achieving a pricing scheme which maximizes the expected profit is difficult because of the discontinuous nature of the profits. We propose a pricing strategy which yields a fixed (say, $\beta$) amount of profit to the charging station. We show that the above  pricing strategy also maximizes the social welfare with the desired level of probability for a suitable choice of $\beta$ (Theorem~\ref{thm:approx}). Hence, such a pricing scheme is also  attractive to the regulators. Further, we show that a suitable choice of $\beta$ can maximize the profit of the charging station for a class of utility functions (Theorem~\ref{thm:aclassutility}). 


Finally, we, empirically provide insights how a trade-off between the profit of the charging station and the social welfare can be achieved for various pricing schemes (Section~\ref{sec:simulation_results}). We also show that how our pricing scheme can increase greater utilization of the resources and result in a lower number of charging spots compared to the existing ones. 

{\em The proofs are deferred to the technical report\cite{tech_ev} owing to the space constraint.}

\textbf{Related Literature}: {\em To the best our knowledge this is the first attempt to consider contract based online pricing for controlling both the energy and deadline of the EVs.} However, other pricing mechanisms to control the charging pattern of EV in a {\em residential} charging  in a day-ahead market are proposed\cite{soltani,oren,javidi,tvt,kar,zou}.  In contrast to  the residential charging, in a commercial or workplace charging station, users do not have control of charging the car at each instance. They need a certain amount of energy within a deadline.  In our proposed mechanism, the charging station selects different prices for different options and the user only needs to select a contract, it does not need to control the charging pattern at each instance. 

Optimal pricing for a day-ahead demand response program have been studied \cite{sen,low2}. However, we need an online pricing mechanism in the EV charging station. Since the charging station is unaware of the future arrivals of the users and has limited renewable energy, determining the optimal pricing in such a setting is more challenging. The menu-based pricing is an online pricing mechanism and can enhance the efficient usage of the resources by controlling the deadline. Unlike in the demand response program, the users also do not need to control the demand for each instance in our menu-based pricing. 

\cite{parkes} proposed  an online VCG auction mechanism. However, in \cite{parkes} the user's payment is determined at the end of the day, and thus users are not sure how much they have to pay beforehand.  Hence, it may not be preferred by the users.  In contrast, in our mechanism the users select one of the contracts  and pay the prescribed price beforehand.   In \cite{tong,qhuang,xu_cdc,yu_allerton} online scheduling algorithms have been proposed for charging EVs.  The main focus of these papers was scheduling, they did not consider the optimal pricing approach for the charging station which we did. Further, unlike in \cite{tong,qhuang,xu_cdc,yu_allerton}, in our menu-based pricing scheme, the charging station can control the energy requirement and the deadline of the users by selecting the prices to the users. Hence, a greater flexibility can be achieved. Additionally, \cite{tong,yu_allerton}  did not guarantee that the energy demand will be fulfilled.  In contrast, in our model once a user opts for an option, the EV charging station always fulfills the request of users.  

In the deadline differentiated pricing \cite{bitar,bitar2,nayyar,Salah2016}, each user's total energy consumption is fixed, however the user can specify the deadline.  On the contrary, in our proposed menu-based pricing mechanism each user can jointly choose {\em any pair} of energy level and deadline. The deadline differentiated pricing is suited for a day-ahead offline setting where an equilibrium price is attained for a specific set of pre-determined decision of the users. However, the users' utilities and thus optimal decision may change in real time and thus, the deadline differentiated pricing may not be suitable for online setting. Our menu-based pricing approach is online, where the price menu is adapted for each arriving user. Further,  \cite{nayyar}  assumed that the price setter knows the utilities of the users. 
\cite{bitar, Salah2016} assumed that the utility functions are strictly concave, and \cite{bitar2} put some restrictions on the utility functions to achieve optimality. However, such assumptions are not necessary for our approach. 

\section{Model}
\begin{figure*}
	\begin{minipage}{0.49\linewidth}
		\includegraphics[trim=1in 0.5in 1in 0in,width=0.99\textwidth]{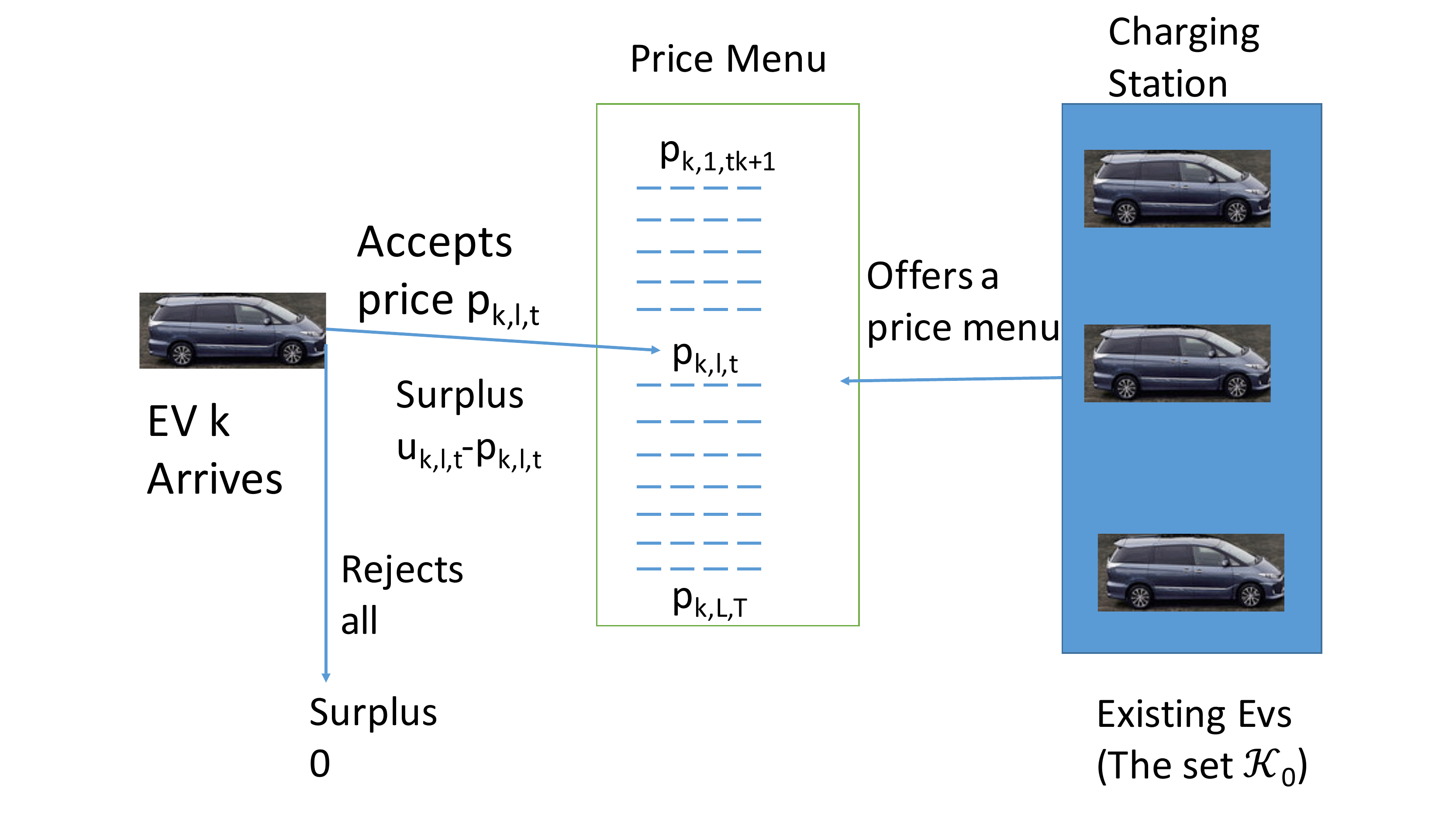}
		\caption{The trading model: Charging station offers a menu of contracts, the arriving user decides either one of them or rejects all.}
		\label{fig:ev_charging}
		\vspace{-0.2in}
	\end{minipage}\hfill
	\begin{minipage}{0.49\linewidth}
		\includegraphics[trim=0.5in 0.7in 1.2in 0in, width=0.99\textwidth]{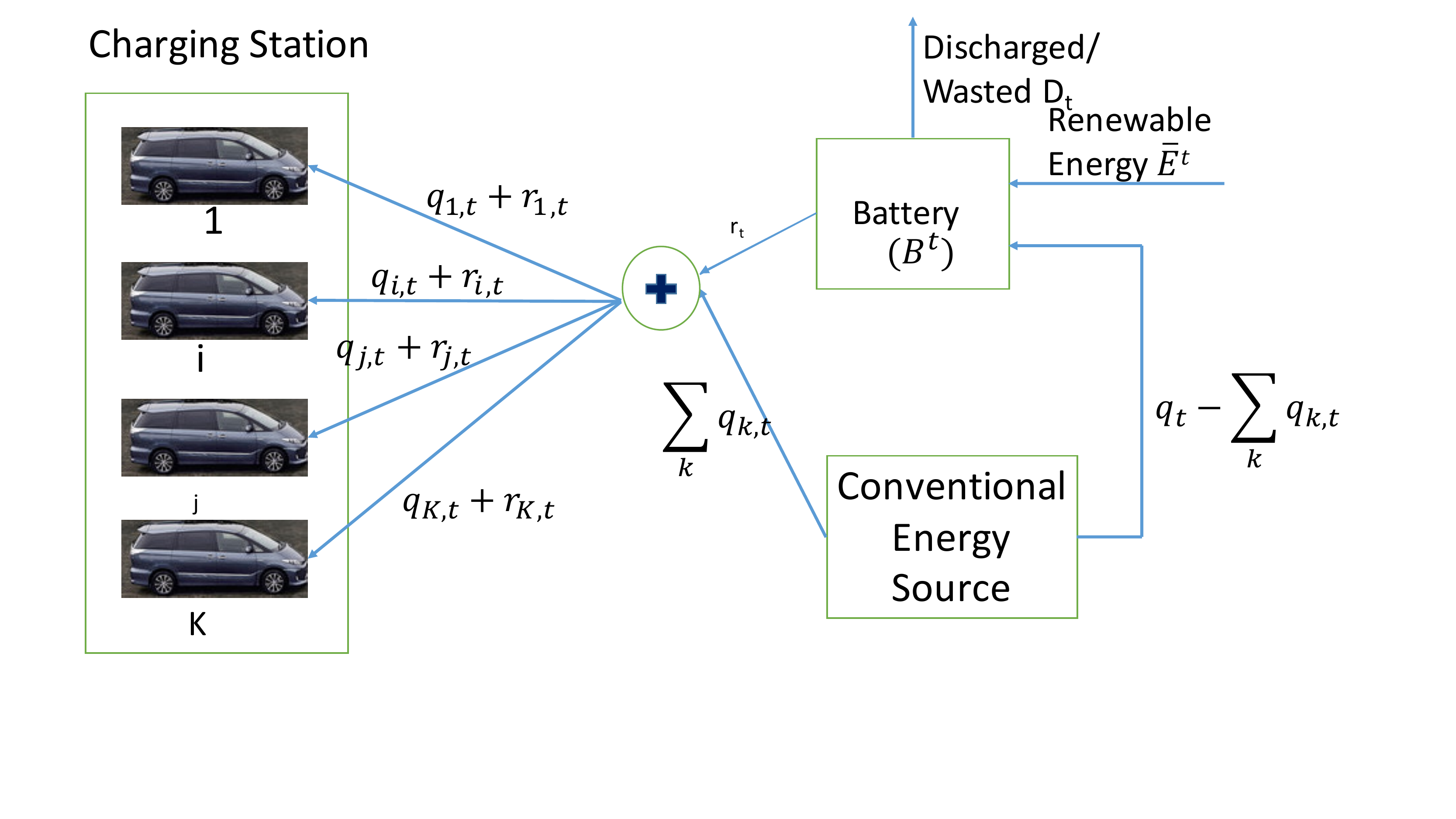}
		\caption{The hybrid energy source with a limited capacity Battery.}
		\label{fig:hybrid_source}
		\vspace{-0.2in}
	\end{minipage}
\end{figure*}

We consider a charging station which wants to select a pricing strategy in order to maximize its payoff over a certain period of time $T$ (e.g. one day).   Suppose that user $k=1,\ldots, K$ arrives at the charging station at time $t_k$. The charging station decides a price menu or a contract $p_{k,l,t}$ to user $k$ for different energy levels $l\in \{1,\ldots, L\}$ and deadline $t\in \{t_{k}+1,\ldots, T\}$ (Fig.~\ref{fig:ev_charging}).\footnote{In Section V, we show that if the support of the distributions of the user's utilities have the same lower end-point, our proposed menu-based price mechanism essentially becomes equivalent to a ime-of-use price mechanism.} It is needless to say that we can discretize the time and energy at any level that one may want\footnote{It is an online mechanism, we can easily extend it for fixed maximum deadline scenario for EVs by extending the time horizon.}, however, the computational cost will increase. User $k$ has to decide $l$ and $t$ based on the menu; if she decided to accept any option on the menu, she has to pay the prescribed price $p_{k,l,t}$. The user can decide not to accept any price too (Fig~\ref{fig:ev_charging}). The EV may not be charged continuously i.e. preemption is allowed. A preempted battery of the EV can be resumed charging from the previous battery level upon preemption.  


\subsection{User's utilities}
If user $k$ selects the price option $p_{k,l,t}$, it will get an utility $u_{k,l,t}$. Hence, its {\em surplus} or payoff will be $u_{k,l,t}-p_{k,l,t}$.  If the user rejects all the options then its utility is $0$ (Fig.~\ref{fig:ev_charging}). We assume that the realized value $u_{k,l,t}$  is drawn from a distribution function of the random variable $U_{k,l,t}$. The random variables $U_{k,l,t}$ need not be independent, in fact, they can be generated from a joint distribution. In practice, there is a correlation of the utilities among different deadlines and charging amount. For example, $U_{k,l_1,t}\geq U_{k,l,t}$ if $l_1>l$ as a higher amount of energy for a fixed deadline should induce higher utility to a user. Similarly, $U_{k,l,t_1}\leq U_{k,l,t_2}$ if $t_1>t_2$ since for a similar level of charge, the user will  prefer the smaller deadline menu as it will give the user more flexibility. On the other hand, a user who  wants to park a long time may not mind a longer deadline. Thus, we do not make any a priori assumptions on the utility functions since they can be different for different users. We assume that the car vacates the charging spot once it exceeds its prescribed deadline\footnote{If the user can not take away its car, the charging spot will be automatically downgraded to a mere parking spot i.e. without any charging facility.}. 
\subsection{The charging Station}
\subsubsection{Hybrid Energy Source}
We assume that the charging station can obtain energies for fulfilling the charging request both from the renewable sources and conventional sources (Fig.~\ref{fig:hybrid_source}).
The charging station can  buy a conventional energy $q_t$ at a price $c_t$ for usage during the interval $[t,t+1)$. We do not assume any specific type of pricing schemes for buying conventional energy, however, we assume that $c_t$ is known. If the real time pricing is used, then we consider $c_t$ as the expected real-time price at time $t$. 

The charging station is also equipped with an energy harvesting device and a storage capacity of $B_{max}$ (Fig.~\ref{fig:hybrid_source}). The harvesting device harvests $\bar{E}^{t}$ amount of energy between $[t,t+1)$. We assume that the marginal cost to harvest renewable energy is $0$. The amount of energy that the charging station uses from the storage for the time $[t,t+1)$ be $r_t$. 

\subsubsection{System Constraints}
The charging station must procure the $l$ amount of energy for user $k$ by time $d_k$ if the user accepts the price menu $p_{k,l,d_k}$.   Let  $q_{k,t}$ be the conventional energy and $r_{k,t}$ be the energy from the storage device used by the charging station to charge user $k$ for time interval $[t,t+1)$. Then, we must have the following constraint
\vspace{-0.07in}
\begin{eqnarray} 
\sum_{t=t_k}^{d_k-1}(r_{k,t}+q_{k,t})\geq l \label{eq:setofaccepted}
\end{eqnarray}
Suppose that  initially, there is a set of users $\mathcal{K}_0$   already present in the charging station  at time $t_k$. Now, the user $i\in \mathcal{K}_0$ has a deadline of $w_i$ and additional demand $N_i$. The charging station must have to satisfy the demand of those users. Thus,  the charging station must also satisfy
\vspace{-0.07in}
\begin{eqnarray}
\sum_{t=t_k}^{w_i-t_k-1}(r_{i,t}+q_{i,t})\geq N_i,\forall i\in \mathcal{K}_0\label{eq:deadline}
\end{eqnarray}
We also assume that the charging station has one kind of charging equipment (either slow charging or fast charging) and there is a maximum rate constraint ($R_{max}$). Hence, 
\vspace{-0.07in}
\begin{align}\label{eq:rate_constraint}
r_{k,t}+q_{k,t}\leq R_{max}, \quad r_{i,t}+q_{i,t}\leq R_{max} \forall i\in \mathcal{K}_0\quad \forall t.
\end{align}
Since the total energy $r_t$ used for charging from the storage device and the amount of conventional energy, $q_t$ bought from the market, thus,
\vspace{-0.07in}
\begin{align}
r_{k,t}+\sum_{i\in \mathcal{K}_0}r_{i,t}=r_t, \quad q_{k,t}+\sum_{i\in \mathcal{K}_0}q_{i,t}\leq q_t\quad \forall t\label{eq:energy}
\end{align}
Note that the charging station may store the unused conventional energy bought from the market i.e. $q_t-q_{k,t}-\sum_{i\in \mathcal{K}_0}q_{i,t}$ is stored in the storage device. The charging station may buy an additional amount of conventional energy at time $t$, if the future prices are higher. 

Let the battery level at time $t+1$ be $B^{t+1}$ and 
$B_{0}$ be the initial battery level. The charging station also wants to keep the battery level at the end of the day as $B_{0}$.  If the final battery level need not match the initial level, our pricing approach can be easily extended to that scenario. If the battery can not hold the excessive energy, then it is wasted. Let use denote the wasted energy for time $[t, t+1)$ be $D_t\geq 0$, \footnote{It is also straightforward to extend our setting when the charging station can sell the excess energy to the grid, i.e. it will sell $D_t$ to the grid.}then
\begin{align}
B^{t+1}= B^{t}+\bar{E}^{t}-r_t+q_t-q_{k,t}-\sum_{i\in \mathcal{K}_0}q_{i,t}-D_t\nonumber\\
0\leq B^{t+1}\leq B_{max},  B^{t_k}=B_{0}\quad B^{T}=B_{0}\label{eq:batterylevel}.
\end{align}
The constraints in (4) and (5) put a bound on the maximum amount of energy can be used for charging.
\section{Problem Formulation}
The profit of the charging station inherently depends on whether the user will accept that menu or not. Hence, before how the charging station will select $p_{k,l,t}$ for user $k$,  we delve into the decision process of the users. 
\subsection{User's decision}\label{sec:userdecision}
A user selects at most one of the price menus in order to maximize its payoff or surplus. We assume that the user is a {\em price taker}. Thus, for a menu of  prices $p_{k,l,t}$, the user $k$ {\em selects} $A_{k,l,t}\in [0,1]$ such that it maximizes the following
\begin{align}\label{eq:utility1}
& \text{maximize} \sum_{l=1}^{L}\sum_{t=t_k+1}^{T}A_{k,l,t}( u_{k,l,t}-p_{k,l,t})\nonumber\\
& \text{subject to } \sum_{l=1}^{L}\sum_{t=t_k+1}^{T}A_{k,l,t}\leq 1
\end{align}
Note from the formulation in (\ref{eq:utility1}) the maximum is achieved  when $A_{k,l,t}=1$  for the contract which maximizes the user $k$'s  payoff (i.e., $\max_{i,j}\{u_{k,i,j}-p_{k,i,j}\}=u_{k,l,t}-p_{k,l,t}$.) and is $0$ otherwise. If such a solution is not unique, any convex combination of these solutions is also optimal since a user can select any of the maximum payoff contracts. We denote the  decision as $A_{k,l,t}(\mathbf{p}_{k})$. Note that the decision whether to accept the menu $p_{k,l,t}$  not only depend on the price $p_{k,l,t}$ but also other price menus i.e. $p_{k,i,j}$ where $i\in \{1,\ldots, L\}$ and $j\in \{t_k+1,\ldots, T\}$ as the user only selects the price menu which is the  most favorable to himself. Note that if the maximum payoff that user gets among all the price menus (or, contracts) is negative, then the user will not charge i.e. $A_{k,l,t}=0$ for all $l$ and $t$. We also assume that if there is a tie between charging and not charging, then the user will decide to charge i.e. if the maximum payoff that user can get is $0$, then the user will decide to charge.\footnote{Our result can be readily extended to the other options, in that case the price strategies given in this paper have to be decreased by an $\epsilon>0$ amount.}

\vspace{-0.1In}
\subsection{Myopic Charging Station}
Since the users arrive for the charging request at any time throughout the day, the charging station does not know the exact arrival times for the future vehicles.  We consider that the charging station is myopic or near-sighted  i.e. it selects its price for user $k$ without considering the future arrival process of the vehicles. However, it will consider the cost incurred to charge the existing EVs. Note that as the number of existing users increases, the marginal cost can increase to fulfill a contract for an arriving user, hence, such a pricing strategy may not maximize the payoff in a long run.  We, later  show that {\em a myopic pricing strategy is optimal in the case the marginal cost of fulfilling a demand of a new user  is independent of the number of existing users.}  

In practice, the charging station often has fixed number of charging spots, thus, the charging station may want to select high prices for user $k$, in order to make the charging spots available for the users who can pay more but only will arrive in future\footnote{The above consideration is left for the future work.}.  However, such a pricing strategy is against the first come first serve basis which is the current norm for charging a vehicle. Our approach considers a fair allocation process, where the charging station serves users based on the first come first serve basis. Later in Section~\ref{sec:simulation_results}, we show that since the charging station can control the time spent by an EV through pricing, our approach results into a lower number of charging spots compared to the existing pricing mechanism.
\subsection{Charging Station's Decisions and cost}\label{sec:chargingstationdecision}
Note that if the user $k$ accepts the menu $(l,d_k)$. Then, the charging station needs to allocate resources among the EVs in order to minimize the total cost of fulfilling the contract. 
First, we introduce some notations which we use throughout. 
\vspace{-0.04in}
\begin{definition}\label{defn:vlt}
	The charging station has to incur the cost $v_{l,d_k}$ for fulfilling the contracts of existing customers and the contract $(l,d_k)$ of the new user $k$, where $v_{l,d_k}$ is the value of the following linear optimization problem:
	\vspace{-0.05in}
	\begin{eqnarray}\label{eq:vlt}
	\mathcal{P}_{l,d_k}:&  \text{min } \sum_{t=t_k}^{T-1}c_tq_t\nonumber\\
	& \text{subject to } (\ref{eq:setofaccepted}),(\ref{eq:deadline}),(\ref{eq:rate_constraint}),(\ref{eq:energy}), (\ref{eq:batterylevel})\nonumber\\
	& \text{var: } r_{k,t}, q_{k,t}, q_t, r_t, D_t\quad \geq 0
	\end{eqnarray}
\end{definition}
Note that our model can also incorporate time varying, strictly increasing convex costs $C_t(\cdot)$\footnote{Our model can also incorporate the scenario where there is a constraint on the maximum value of conventional energy can be bought from the market at a given time.}. Since $\mathcal{P}_{l,d_k}$ is a linear optimization problem, it is easy to compute $v_{l,d_k}$. Also, note that if the above problem is infeasible for some $l$ and $t$, then we consider $v_{l,t}$ as $\infty$. We assume that the prediction $\bar{E}^{t}$ is perfect for all future times and known to the charging station\footnote{ Menu-based pricing approach can be extended to the setting where the estimated generation does not match the exact amount. First, we can consider a conservative approach where $\bar{E}^t$ can be treated as the worst possible renewable energy generation. As a second approach, we can  accumulate various possible  scenarios of the renewable energy generations, and try to find the cost to fulfill a contract for each such scenario. For example, if there are $M$ number of possible instances of the renewable energy generation amount in future. Then, we can find the optimal cost for each such instance of renewable energy generation $\bar{E}^{m,t}$ where $m\in \{1,\ldots, M\}$ instead of $\bar{E}^t$ . We then can compute the average (or, the weighted average, if some instance has greater probability) of the optimal costs, and that cost can be taken as the cost of fulfilling a certain contract.}.
\vspace{-0.04in}
\begin{definition}\label{defn:v-k}
	Let $v_{-k}$ be the amount that the charging station has to incur to satisfy the requirements of the existing EVs if the new user does not opt for any of the price menus. 
\end{definition}
If user $k$ does not accept any price menu, then the charging station still needs to satisfy the demand of existing users i.e. the charging station must solve the problem $\mathcal{P}_{l,t}$ with $q_{k,t}=r_{k,t}=0$. $v_{-k}$ is the value of that optimization problem.   Thus, from Definitions~\ref{defn:vlt} and \ref{defn:v-k} we can visualize $v_{l,d_k}-v_{-k}$ as the additional cost or marginal cost to the charging station when the user $k$ accepts the price menu $p_{k,l,d_k}$. It is easy to discern that {\em $v_{l,d_k}-v_{-k}$ is non-negative for any $d_k$ and $l$.} 



\subsection{Profit of the charging station}\label{sec:profit}
Now, we discuss the profit  of the charging station based on its pricing strategies. Note that  if all the spots are occupied then, the charging station can not accommodate a new user.  Thus, we consider the scenario where  a charging spot is available. 
\subsubsection{Pricing with Perfect Foresight}
First, we consider the scenario where the charging station has a perfect foresight of the utility of the user i.e. the charging station is clairvoyant and has a perfect knowledge about the user's utility.   Note that if the user $k$ selects the price menu $p_{k,l,t}$, then the charging station has to pay  $v_{l,t}$ amount (Definition~\ref{defn:vlt}). Thus, the charging station has to pay additional amount $v_{l,t}-v_{-k}$ (Definition~\ref{defn:v-k}) when the user selects the price menu $p_{k,l,t}$.  Thus, the profit of the charging station is 
\vspace{-0.05in}
\begin{eqnarray}\label{eq:profit_maximum}
\sum_{l=1}^{L}\sum_{t=t_k+1}^{T}(p_{k,l,t}-v_{k,l,t}+v_{-k})A_{k,l,t}(\mathbf{p}_{k})
\end{eqnarray}
Note that here the charging station selects a price $p_{k,l,t}$ to maximize its own profit. $A_{k,l,t}>0$ only if $p_{k,l,t}$ gives the highest payoff for the user $k$ as discussed in Section~\ref{sec:userdecision}.
\subsubsection{Prediction Based Pricing}
In practice the charging station may not know the exact realization of the utility function of the users. Thus, it can only use predictions of the utility function in order to select the price menu. Here, we consider such a scenario where the charging station does not know the exact utilities of the user. 

We assume that the charging station knows the statistic of the user's utility. Let $R_{k,l,t}$ be the event that the price menu $p_{k,l,t}$ is selected, hence, the profit of the charging station is
\vspace{-0.05in}
\begin{align}\label{eq:profitmax_expec}
& \sum_{l=1}^{L}\sum_{t=t_k+1}^{T}\mathbb{E}[(p_{k,l,t}-v_{l,t}+v_{-k})\mathbbm{1}_{R_{k,l,t}}]\nonumber\\
& =\sum_{l=1}^{L}\sum_{t=t_k+1}^{T}(p_{k,l,t}-v_{l,t}+v_{-k})\Pr(R_{k,l,t})
\end{align}
The indicator variable $R_{k,l,t}$ in (9) denotes the event that the contract $(l,t)$ is being chosen by the user $k$.
The expectation is taken over the joint distribution of $U_{k,l,t}$. {\em The expected profit maximization problem for the charging station is to maximize the above objective over} $p_{k,l,t}$. 
We assume that the utilities are distributed according to some continuous distribution\footnote{However, it can be easily extended to discontinuous distribution case.}. Hence, $\Pr(R_{l,t})$ is given by
\vspace{-0.05in}
\begin{align}
\Pr(R_{k,l,t})= \Pr(U_{k,l,t}-p_{k,l,t}\geq (\max_{i,j}\{ U_{k,i,j}-p_{i,j}\})^{+})\nonumber
\end{align}
Thus, $\Pr(R_{k,l,t})$ not only depends on $p_{k,l,t}$, but also prices for other menus. 
\subsection{Objectives}
We consider that the charging station decides the price menus in order to fulfill one of the two objectives (or, both)--i) Social Welfare Maximization and ii) its profit maximization.

\subsubsection{Social Welfare}
The social welfare is the sum of  user surplus and the profit of the charging station. As discussed in Section~\ref{sec:userdecision} for a certain realized values $u_{k,l,t}$ if the user $k$ selects the price menu $p_{k,l,t}$, then its surplus is $u_{k,l,t}-p_{k,l,t}$, otherwise, it is $0$.

As discussed in Section~\ref{sec:profit} the profit of the charging station is $p_{k,l,t}-v_{l,t}+v_k$ for a given price $p_{k,l,t}$  if the user selects the menu, otherwise it is $0$.  Hence, the social welfare maximization problem is to select the price menu $p_{k,l,t}$ which will maximize the following
\begin{align}\label{eq:socialwelfare_maxproblem}
\mathcal{P}_{perfect}: & \text{maximize} \sum_{l=1}^{L}\sum_{t=t_k+1}^{T}(u_{k,l,t}-v_{l,t}+v_{-k})A_{k,l,t}(\mathbf{p}_k)\nonumber\\
& \text{var }: p_{k,l,t}\geq 0.
\end{align}
Recall that in order to find $v_{l,t}$ we have to solve $\mathcal{P}_{l,t}$(cf.~(\ref{eq:vlt})) which is a constrained optimization problem. 
Since EV is expected to increase the social value such as providing a cleaner environment, and higher energy efficiency, hence, it is  important for a regulator (e.g. FERC) whether there exists a pricing strategy which maximizes the social welfare of the system.  If the charging station is operated by the regulator or some government agency, then the main objective is indeed maximizing the social welfare or user surplus is maintained.

{\em Ex-ante and Ex-post Maximization}: When the charging station is unaware of the utilities of the users, then two options are considered--i) decides a price and hopes that it will maximize the social welfare for the realized values of utilities ({\em ex-post} maximization), or ii) decides a price and hopes that it will maximize  the social welfare in an expected sense ({\em ex-ante} maximization). 
Thus,  the {\em ex-ante} maximization  does not guarantee that the social welfare will be maximized for every realization of the random variables $U_{k,l,t}$. However, in the {\em ex-post} maximization, the social welfare is maximized for each realization of the random variables. Thus, {\em ex-post} maximization is a stronger concept of maximization (and thus, more desirable) and it is not necessary that there exist pricing strategies which maximize the ex-post  social welfare. However, we show that in our setting there exist pricing strategies which maximize the ex-post social welfare. Note that {\em ex-post} social welfare maximization is the same as (\ref{eq:socialwelfare_maxproblem}).

\subsubsection{Profit Maximization}
Social welfare maximization does not guarantee that the charging station may get a positive profit. It is important for the wide-scale deployment of the charging stations, the charging station must have some profit. Further, if the charging station is operated by a private entity its objective is indeed to maximize the profit. 

When the charging station is clairvoyant, then the charging station wants to maximize the profit given in (\ref{eq:profit_maximum}) by selecting $p_{k,l,t}$. On the other hand when the charging station does not know the user's utility, then it wants to maximize the expected payoff  given in (\ref{eq:profitmax_expec}) by selecting $p_{k,l,t}$.

\subsubsection{Separation Problem}
Note that in order to select optimal $p_{k,l,t}$, the charging station has to obtain $v_{l,t}$  and $v_{-k}$ (Definitions~\ref{defn:vlt} \& \ref{defn:v-k}). However,  $v_{l,t}$ and $v_{-k}$  do not depend on $p_{k,l,t}$. Hence, we can separate the problem--first the charging station finds $v_{l,t}$ and $v_k$, and then it will select $p_{k,l,t}$ to fulfill the objective. {\em We now focus on finding optimal} $p_{k,l,t}$.
\section{Results: Social Welfare Maximization}
First, we state the optimal values of the social welfare for any given realization of the user's utilites. Next, we state a pricing strategy which attains the above optimal value. 

Note that if $u_{k,l,t}-v_{l,t}+v_{-k}<0$ for each $l$ and $t$, then the social welfare is maximized when the user $k$ does not charge or equivalently, the price $p_{k,l,t}$ is very high  for each $l$ and $t$. In this case, the optimal value of social welfare is $0$.

On the other hand if $u_{k,l,t}-v_{l,t}\geq -v_{-k}$ for  some $l$ and $t$, then the  social welfare is maximized when  the user $k$ charges its car. If the user accepts the price menu $p_{k,l,t}$, then the social welfare is $u_{k,l,t}-v_{l,t}+v_{-k}$. Thus, the maximum social welfare in the above scenario is $\max_{l,t}(u_{k,l,t}-v_{l,t}+v_{-k})$. Hence-
\vspace{-0.03in}
\begin{theorem}\label{thm:max_socialwelfare}
	The maximum value of social welfare is $\max\{\max_{l,t}(u_{k,l,t}-v_{l,t}+v_{-k}),0\}$.
\end{theorem}
Note that even though the maximum value of social welfare is unique (as in Theorem~\ref{thm:max_socialwelfare}), the optimal pricing strategy is not unique. In the following, we give one possible pricing strategy that achieves the optimal social welfare.
\vspace{-0.04in}
\begin{theorem}\label{thm:pricestrategy}
	The pricing strategy $p_{k,l,t}=v_{l,t}-v_{-k}$ maximizes the social welfare. 
\end{theorem}
Note that in the pricing strategy, the charging station does not need to know the utility of the users. It optimizes the social welfare for each possible realization of the utility functions. Hence, we obtain 
\vspace{-0.04in}
\begin{corollary}\label{thm:price_expected}
	The pricing strategy $p_{k,l,t}=v_{l,t}-v_{-k}$ maximizes the ex-post social welfare. 
\end{corollary}
Though the pricing strategy maximizes the social welfare, the above pricing strategy does not provide any positive profit to the charging station. Thus, the charging station may not prefer this pricing strategy as it will not have any incentive to provide the charging spots. 

Also note that the pricing strategy also maximizes the social welfare in the long run when the additional cost of fulfilling a contract (i.e. $v_{l,t}-v_{-k}$) does not depend on the existing users in the charging station. 
The condition that $v_{l,t}-v_{-k}$ is independent of the existing EVs in the charging station is satisfied if either all demand can be fulfilled using renewable energy or there is no renewable energy generation and the conventional energy is bought at a flat rate. Hence, in the {\em two above extreme cases, the myopic pricing strategy is also optimal in the long run}. 

\section{Profit Maximization of the charging station}
\subsection{Charging station with perfect foresight}
We now provide a price strategy which maximizes the profit of the charging station and also the social welfare when {\em the charging station is clairvoyant}. Recall that the profit of the charging station is given by (\ref{eq:profit_maximum}). 

First, we introduce a notation.
\vspace{-0.05in}
\begin{definition}\label{defn:optimal}
	Let $(l^{*},t^{*})=$argmax$_{l,t}\{u_{k,l,t}-v_{l,t}\}$.
\end{definition}
\vspace{-0.03in}
\begin{theorem}\label{thm:profit_max}
	Set $p_{k,l,t}=v_{l,t}-v_{-k}+(u_{k,l^{*},t^{*}}-v_{l^{*},t^{*}}+v_{-k})^{+}$ where $(l^{*},t^{*})$ is defined in Definition~\ref{defn:optimal}. Such a pricing strategy maximizes the profit as well as the social welfare.
\end{theorem}
The above pricing strategy is an example of {\em value-based} pricing strategy where prices are set depending on the valuation or the utility of the users \cite{Hinterhuber}. The user's utility dependent pricing strategy is also proposed in smart grids in some recent papers \cite{kar,personalized_pricing}. In contrast, the price strategy stated in Theorem~\ref{thm:pricestrategy} is an example of {\em cost-based} pricing strategy where the prices only depend on the costs. If the utilities of the users are same, the pricing strategy becomes similar to a time-dependent pricing scheme, which is prevalent in practice. 

In the value-based pricing strategy, the user surplus decreases, in fact it is\footnote{If the user is reluctant to charge if it does not get a positive payoff, then, we can reduce the price by $\epsilon>0$ amount. In that case, it will be $(1-\epsilon)$ optimal profit-maximizing strategy.} $0$ in our case.  Thus all the user surplus is transferred as the profit of the charging station. Hence, when the charging station is clairvoyant, then the pricing strategy which maximizes the profit of the charging station and it does not entail any positive {\em user surplus}. 

{\em Note that there can be other  pricing strategies which simultaneously maximize the social welfare and the profit.} For example, if $p_{k,l,t}$ is $\infty$ for all $(l,t)\neq (l^{*},t^{*})$ and $p_{k,l^{*},t^{*}}=v_{l^{*},t^{*}}-v_{-k}+(u_{k,l^{*},t^{*}}-v_{l^{*},t^{*}}+v_{-k})^{+}$, then it also maximizes the profit of the charging station. Thus, in this scenario, it can give only one possible contract to the EVs.

Though the joint  profit maximizing and social welfare pricing strategy may not be unique, the profit of the charging station is the unique and  is given by 
\begin{align}
(u_{k,l^{*},t^{*}}-v_{l^{*},t^{*}}+v_{-k})^{+}
\end{align}

\subsection{Prediction based pricing}\label{sec:price_uncertainty}

\subsubsection{Maximum Profit under ex-post social welfare maximization}
Note from Theorem~\ref{thm:profit_max} that the profit maximization pricing strategy which maximizes the social welfare requires that the charging station has the complete information of the utilities of the users. Hence, such a pricing strategy can not be implemented when the charging station does not know the exact utilities of the users. Note from (\ref{eq:profitmax_expec})  that the profit maximization is a difficult problem as the user will select one menu inherently depends on the prices selected for other menus. For example, if the price selected for a particular contract is high, the user will be reluctant to take that compared to a lower price one. The profit is a discontinuous function of the prices and thus, the problem may not be convex even when the marginal distribution of the utilities are concave.

We have already seen (Corollary 1)  that a pricing strategy which can maximize the ex-post social welfare, however, it does not give any positive profit. 
We now  show that there exists a pricing strategy which may provide better profit to the charging station while maximizing the ex-post social welfare. First, we introduce a notation which we use throughout.
\begin{definition}\label{defn:lower_endpoint}
	Let $L_{k,l,t}$ be the lowest end-point of the marginal distribution of the utility $U_{k,l,t}$.
\end{definition}
\vspace{-0.04in}
\begin{theorem}\label{thm:profitmax_uncertainty}
	Consider the pricing strategy: 
	\begin{align}\label{eq:pricesocialmax}
	p_{k,l,t}=v_{l,t}-v_{-k}+(\max_{i,j}\{L_{k,i,j}-v_{i,j}+v_{-k}\})^{+}.
	\end{align} The pricing strategy maximizes the ex-post social welfare.

	The profit  is $(\max_{i,j}\{L_{k,i,j}-v_{i,j}+v_{-k}\})^{+}$.
\end{theorem}
The pricing strategy maximizes the ex-post social welfare similar to Corollary~\ref{thm:price_expected}. This is also the {\em maximum possible profit that the charging station can have under the condition that it  maximizes the ex-post social welfare with probability $1$}. However, it may not maximize the expected profit of the charging station or in other words, the pricing strategy which maximizes the expected profit needs not maximize the ex-post social welfare.  Hence, unlike in the scenario where the charging station is clairvoyant (Theorem~\ref{thm:profit_max}) {\em there may not exist a profit maximization strategy which is also a social welfare maximizer when the charging station is unaware of the utilities}. Note that the user surplus is {\em not} $0$, hence,  uncertainty regarding the user's utility functions is required for a positive consumer's surplus.  

Also, note  the similarity with Theorem~\ref{thm:profit_max}. If the user knows the utility, then $L_{k,l,t}=u_{k,l,t}$ as there is no uncertainty and we get back the pricing strategy stated in Theorem~\ref{thm:profit_max}. 

Note that if $\max_{l,t}(L_{k,l,t}-v_{l,t}+v_{-k})>0$, then such a pricing strategy gives  a positive profit to the charging station. If the charging station has large storage or large renewable energy harvesting devices, then, the cost $v_{l,t}-v_{-k}$ will be lower and thus, the charging station can get a higher profit.  It also increases the user  surplus, as the price set by the charging station decreases. Thus, the impact of higher degrees of renewable energy  integration for the charging station increases both the profit of the charging station and the user surplus. The above illustrates the importance of the storage and harvesting energy devices in the charging station. The regulator (e.g. FERC) can also provide incentives to the charging station to set up those devices as the pricing strategy increases profit to the charging station as well as  the ex-post social welfare. 

In the extreme, when $v_{l,t}=0$ for all $l$ and $t$, then the profit of the charging station becomes maximum.  However,  further decreasing $v_{l,t}$ will not have any effect on the profit of the charging station as well as the user surplus, thus, it also shows the investment that the charging station needs to make for storage and renewable energy harvesting devices.

Also note that the users which have higher utilities i.e. $L_{k,l,t}$ is higher, it will give more profits to the charging station. 

The charging station needs to know the lowest end-points of the support set of the utilities unlike in Corollary~\ref{thm:price_expected}. However, the charging station does not need to know the exact distribution functions of the utilities. The lowest end-point can be easily obtained from the historical data. For example, $L_{k,l,t}$ may be computed by the lowest possible price that the user accepts for the energy level $l$ and the deadline $t$. 


\subsubsection{Guaranteed positive profit to the Charging station}
In Theorem~\ref{thm:profitmax_uncertainty} the charging station only has a positive profit if $\max_{l,t}\{L_{k,l,t}-v_{l,t}+v_k\}>0$. In the case, the above condition is not satisfied, then the charging station's profit will be $0$. Naturally, the question is whether there exists a pricing strategy which gives a guaranteed positive profit to the charging station without decreasing  the social welfare much. In the following we provide such a pricing strategy. 

First note that by the continuity of the joint distribution function we have the following
\vspace{-0.04in}
\begin{lemma}\label{lm:delta}
	Let for each $\epsilon>0$, there exists a $\delta>0$ such that 
	\begin{align}
	& \Pr(\max_{l,t}\{U_{k,l,t}-v_{l,t}+v_{-k}\}\geq 0)\nonumber\\
	& \leq \epsilon+\Pr(\max_{l,t}\{U_{k,l,t}-v_{l,t}+v_{-k}-\delta\}\geq 0)
	\end{align} 
\end{lemma}
\vspace{-0.04in}
\begin{theorem}\label{thm:approx}
	Fix an $\epsilon>0$. Now, consider the pricing strategy
	\vspace{-0.12in}
	\begin{align}\label{eq:approx_price}
	p_{k,l,t}=v_{l,t}-v_{-k}+\delta(\epsilon)
	\end{align}
	where $\delta(\epsilon)$ is the $\delta$ which satisfies the Lemma~\ref{lm:delta}.
	
	Then such a pricing strategy maximizes  the ex-post social welfare with probability $1-\epsilon$. 
\end{theorem}
{\em Outline of proof}: First, note that adding a constant does not change the optimal solution. Hence, if  $(l^{*},t^{*})=\text{argmax}_{l,t}(u_{k,l,t}-v_{l,t}+v_{-k})$, then $(l,^{*},t^{*}) $ is also optimal for price strategy in (\ref{eq:approx_price}).  The  rest of the proof follows from Lemma~\ref{lm:delta}. Lemma~\ref{lm:delta} entails that there exists some $\delta>0$ such that it will ensure that the price strategy is off from the social welfare maximizer pricing strategy by at most $1-\epsilon$ in probability. \qed

Note that the pricing strategy stated in (\ref{eq:approx_price}) gives a positive profit of $\delta(\epsilon)$ amount irrespective of the menu selected by the user. Note that {\em the assumption of a continuous distribution is key}. If the distributions are discrete, then $\delta(\epsilon)$ may be $0$. Hence, the charging station may  get zero profit.


The expected profit of the charging station for the above pricing strategy can be readily obtained--
\vspace{-0.04in}
\begin{theorem}\label{thm:expected_payoff}
	Suppose that $p_{k,l,t}=v_{l,t}-v_{-k}+\beta$, then the expected profit of the charging station is $\beta\max_{l,t}\{\Pr(U_{k,l,t}\geq v_{l,t}-v_{-k}+\beta)\}$.
\end{theorem}
{\em Outline of the Proof}: Note that if a user selects any of the contracts, then the charging station's profit is $\beta$. Hence, the charging station's expected profit is $\beta$ times the probability that at least one of the contracts will be accepted.\qed

The regulator such as FERC can select  $\beta$ judiciously to trade off between the profit of the charging station and the user surplus. We empirically study the effect of $\beta$ in Section~\ref{sec:simulation_results}.

Now, we provide an example where the above pricing strategy can also be a profit maximizing for a suitable choice of $\beta$.  First, we introduce a notation
\vspace{-0.06in}
\begin{definition}\label{defn:alpha}
	Let $\zeta=\max\{\gamma| \gamma=\text{argmax}_{\beta\geq 0}\beta\{\max_{i,j}\Pr(U_{k,i,j}>\beta+v_{i,j}-v_{-k})\}\}$. 
\end{definition}
Note that since $U_{k,l,t}$ is bounded and the probability distribution is continuous, thus, $\zeta$ exists. Note from Theorem~\ref{thm:expected_payoff} that $\zeta$ corresponds to the $\beta$ for which the charging station can get the maximum possible profit when the prices are of the form $p_{k,l,t}=v_{l,t}-v_{-k}+\beta$.  

Now consider a class of widely seen utility functions
\vspace{-0.04in}
\begin{assumption}\label{assum:utility}
	Suppose that the utility function $U_{k,l,t}=(Y_{k,l,t}+X_k)^{+}$ for all $l$ \& $t$; $Y_{k,l,t}$ is a constant and known to the charging station, however, $X_k$ is a random variable and whose realized value is not known to the charging station.
\end{assumption}
In the above class of utility function, the uncertainty is only regarding the realized value of the random variable $X_k$. Note that $X_k$ is independent of $l$ and $t$, hence,$X_k$ is considered to be an additive white noise.

It is important to note that we do not put any assumption whether {\em $X_k$ should be drawn from a continuous or discrete distribution.} However, if the distribution is discrete, we need the condition that $\zeta$ must exist. 
\vspace{-0.07in}
\begin{theorem}\label{thm:aclassutility}
	Consider the pricing strategy 
	$p_{k,l,t}=v_{l,t}-v_{-k}+\zeta$;
	where $\zeta$ is defined in Definition~\ref{defn:alpha},
	
	The pricing strategy maximizes the expected profit of the charging station (given in (\ref{eq:profitmax_expec})) when the utility functions are of the form given in Assumption~\ref{assum:utility}. 
\end{theorem}
{\em Remark}: The above result is surprising. It shows that a simple pricing mechanism such as fixed profit can maximize the expected payoff for a large class of utility functions. However, if the utilities do not satisfy Assumption~\ref{assum:utility} then, the above pricing strategy may not be optimal.
\subsection{The pricing algorithm}
\begin{enumerate}
	\item User $k$ comes at time $t_k$.
	\item The charging station solves the linear programming problem $\mathcal{P}_{l,t}$ (eq. (7)) and finds the additional cost $v_{l,t}-v_{-k}$ for fulfilling the contract $(l,t)$ for user $k$ for each $l$ and $t$. 
	\item The charging station selects the price $p_{k,l,t}=v_{l,t}-v_{-k}+\max_{i,j}\{L_{k,i,j}-v_{i,j}+v_{-k}\}^{+}+\beta.$ where $\beta\geq 0$. 
	\item The user selects the contract which maximizes its payoff (eq.(6)). 
\end{enumerate}
Note that when $\beta=0$ gives the worst possible payoff to the charging station as discussed before. The charging station needs to solve the linear programming problem $\mathcal{P}_{l,t}$. The linear programing problem can be efficiently solved by many solvers such as MOSEK, CPLEX, Simplex, CVX, and Linprog tool of MATLAB.

\section{Simulation Results}\label{sec:simulation_results}
We numerically study and compare various pricing strategies presented in this paper. We evaluate the profit of the charging station and the user's surplus  achieved in those pricing strategies. We also show that our mechanism requires  less charging spots   compared to our nearest pricing model.
\vspace{-0.2in}
\subsection{Simulation Setup}
Similar to \cite{quadratic}, the user's utility  for energy $x$  is taken to be of the form $\min\{-ax^2+bx, \dfrac{b^2}{4a}\}$.
Thus, the user's utility is a strictly increasing and concave function in the amount energy consumed $x$. The quadratic utility functions for EV charging have also been considered in \cite{low}. Note that the user's desired level of charging is $b/(2a)$. {\em We assume that $b/(2a)$ is a random variable}. \cite{gov} shows that in a commercial charging station, the average amount of energy consumed per EV is $6.9$kWh with standard deviation $4.9$kWh. We thus consider that $b/(2a)$ is a truncated Gaussian random variable with mean $6.9$kWh and standard deviation $4.9$kWh in the interval $[2, 20]$. We assume $a$ is a uniform random variable in the interval $[1/20, 1/8]$. 

From \cite{gov}, the deadline or the time spent by an electric vehicle in a commercial charging is distributed with an exponential distribution with mean $2.5$ hours. Thus, we also consider the preferred deadline ($T_{pref}$) of the user  to be an exponentially distributed random variable with mean $2.5$. The user strictly prefers a lower deadline. Hence, we assume that the utility is a convex decreasing function of the deadline \cite{tong}. The utility of the user after the preferred deadline is considered to be $0$. Hence, the user's utility is
\begin{align}\label{eq:simulation}
& U_{k,l,t}=\min\{-al^2+lb, b^2/4a\}\times\nonumber\\
& (\exp(T_{pref}-t-t_k)-1)^{+}/(\exp(T_{pref}-t_k)-1)
\end{align}
The arrival process of electric vehicles is considered to be a Poisson arrival process. However, the arrival rates vary over time.  For example, during the peak-hours (8 am to 5pm) the arrival rate is higher compared to the off-peak hours. We, thus, consider a non-homogeneous Poisson process with the arrival rate is $15$  ($5$, resp.) vehicles per hour during the peak period (off-peak period, resp.). We also assume that the maximum charging rate $R_{max}$ is $3.3$ Kw.

We assume that the renewable energy is harvested according to a truncated Gaussian distribution with mean $2$ and variance $2$ per hour. The storage unit is assumed to be of capacity $20$kW-h. Initial battery level is assumed to be $0$ i.e. it is fully discharged.  The prices for the conventional energy is assumed to be governed by Time-of-Use (ToU) time scale. Thus,  the cost of buying the conventional energy varies over time. 
\subsection{Results}
We consider the scenario where the charging station is unaware of the exact utilities of the users, however, it knows the distribution function. We consider the pricing strategy that we have introduced in Section~\ref{sec:price_uncertainty}--
\begin{eqnarray}
p_{k,l,t}=v_{l,t}-v_{-k}+\max_{i,j}\{L_{k,i,j}-v_{i,j}+v_{-k}\}^{+}+\beta.\nonumber
\end{eqnarray}
Recall from Definition~\ref{defn:lower_endpoint} that $L_{k,l,t}$ is the lowest end-point of the utility $U_{k,l,t}.$ We study the impact of $\beta$. 
\subsubsection{Effect on Percentage of the users admitted}
Fig.~\ref{fig:betavsadmittedusers} shows that as $\beta$ increases the number of admitted users decreases. However, the decrement is slow initially. As $\beta$ becomes  larger than a threshold, the price selected to the users becomes very large, and thus, a fewer number of EVs are admitted.
\subsubsection{Effect of $\beta$ on User's Surplus and Profit of the charging station} 
\begin{figure*}
	\begin{minipage}{0.24\linewidth}
		\includegraphics[trim=0in 0in 1in 0in,width=\textwidth]{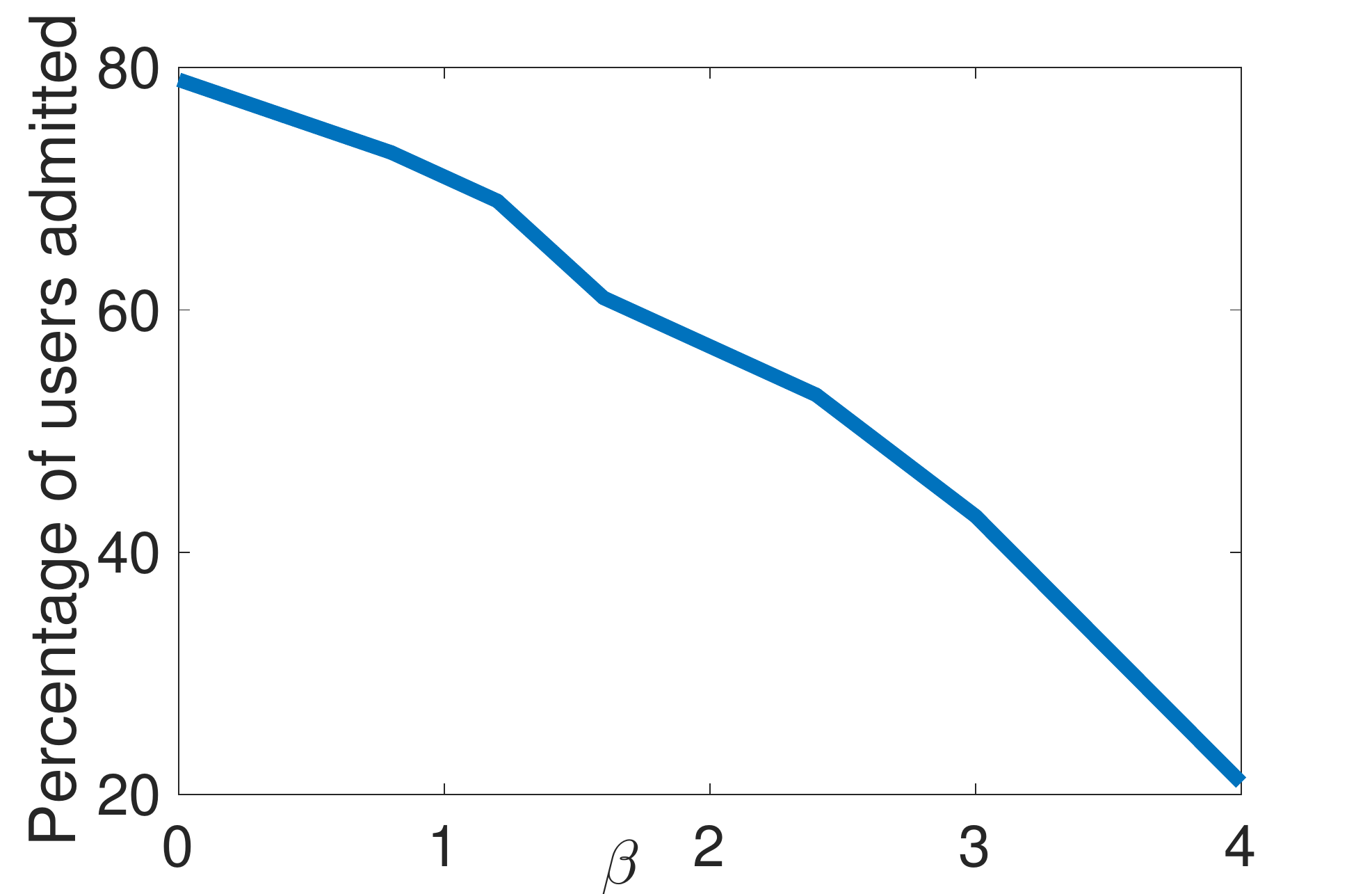}
		\vspace{-0.3in}
		\caption{Variation of the percentage of EVs admitted with $\beta$.}
		\label{fig:betavsadmittedusers}
		\vspace{-0.2in}
	\end{minipage}\hfill
	\begin{minipage}{0.24\linewidth}
		\includegraphics[trim=0in 0in .5in 0in,width=\textwidth]{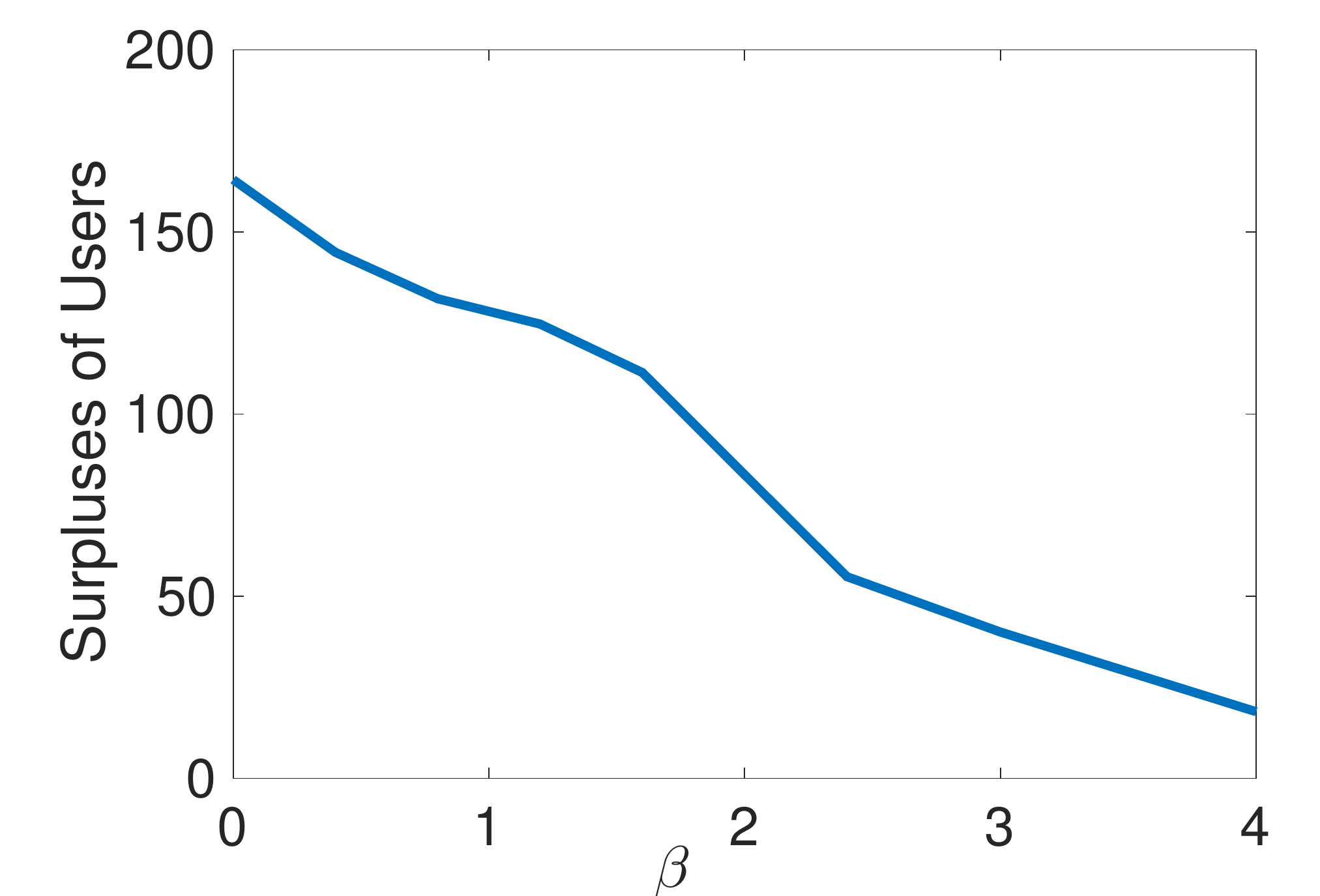}
		\vspace{-0.3in}
		\caption{Variation of the total users' surpluses with $\beta$.}
		\label{fig:consumer_surplus}
		\vspace{-0.2in}
	\end{minipage}
	\begin{minipage}{0.24\linewidth}
		\includegraphics[trim=0in 0in .5in 0in,width=\textwidth]{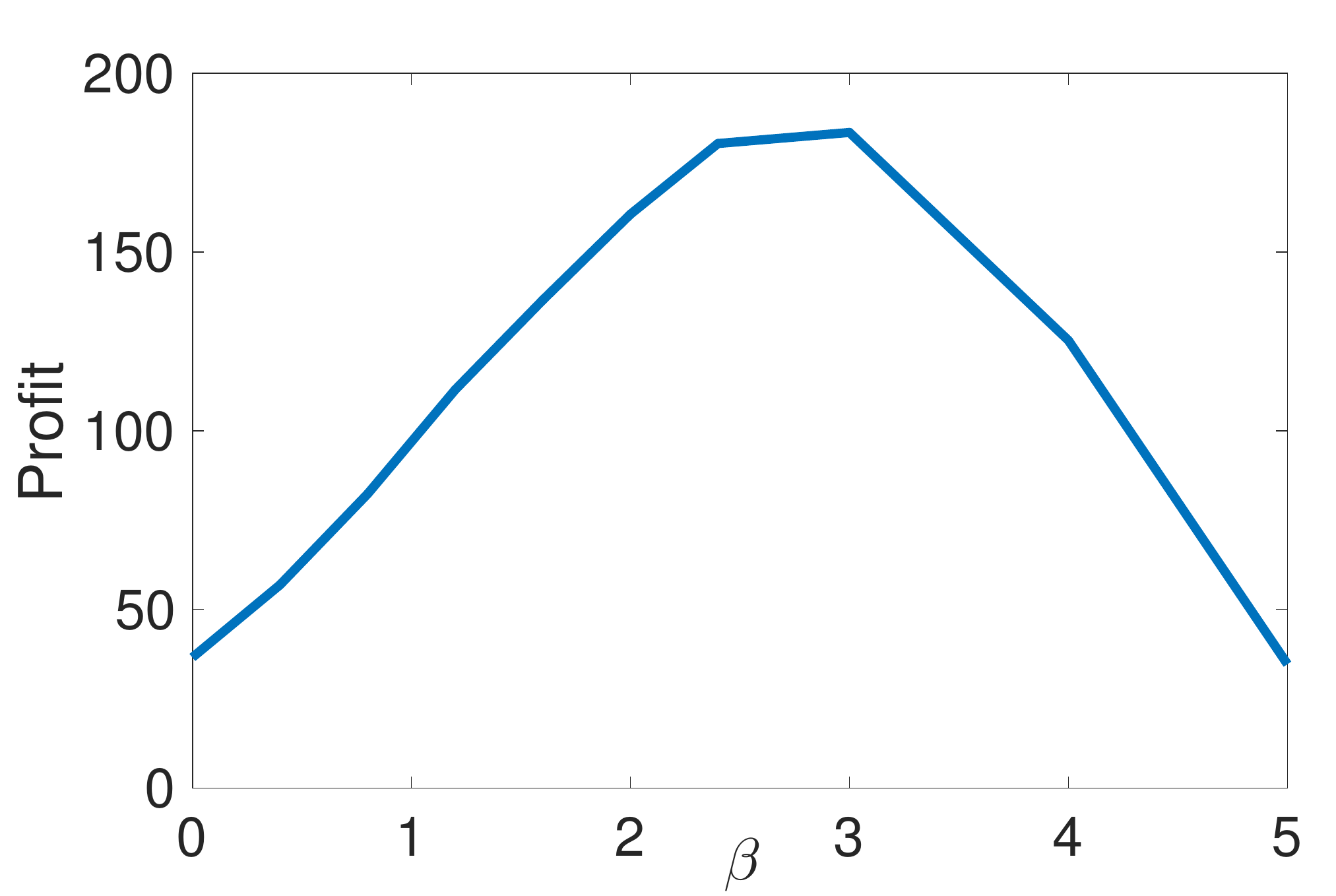}
		\vspace{-0.3in}
		\caption{Variation of the profit of the charging station with $\beta$.}
		\label{fig:profit}
		\vspace{-0.2in}
	\end{minipage}\hfill
	\begin{minipage}{0.24\linewidth}
		\includegraphics[trim=0in 0in .7in 0in,width=\textwidth]{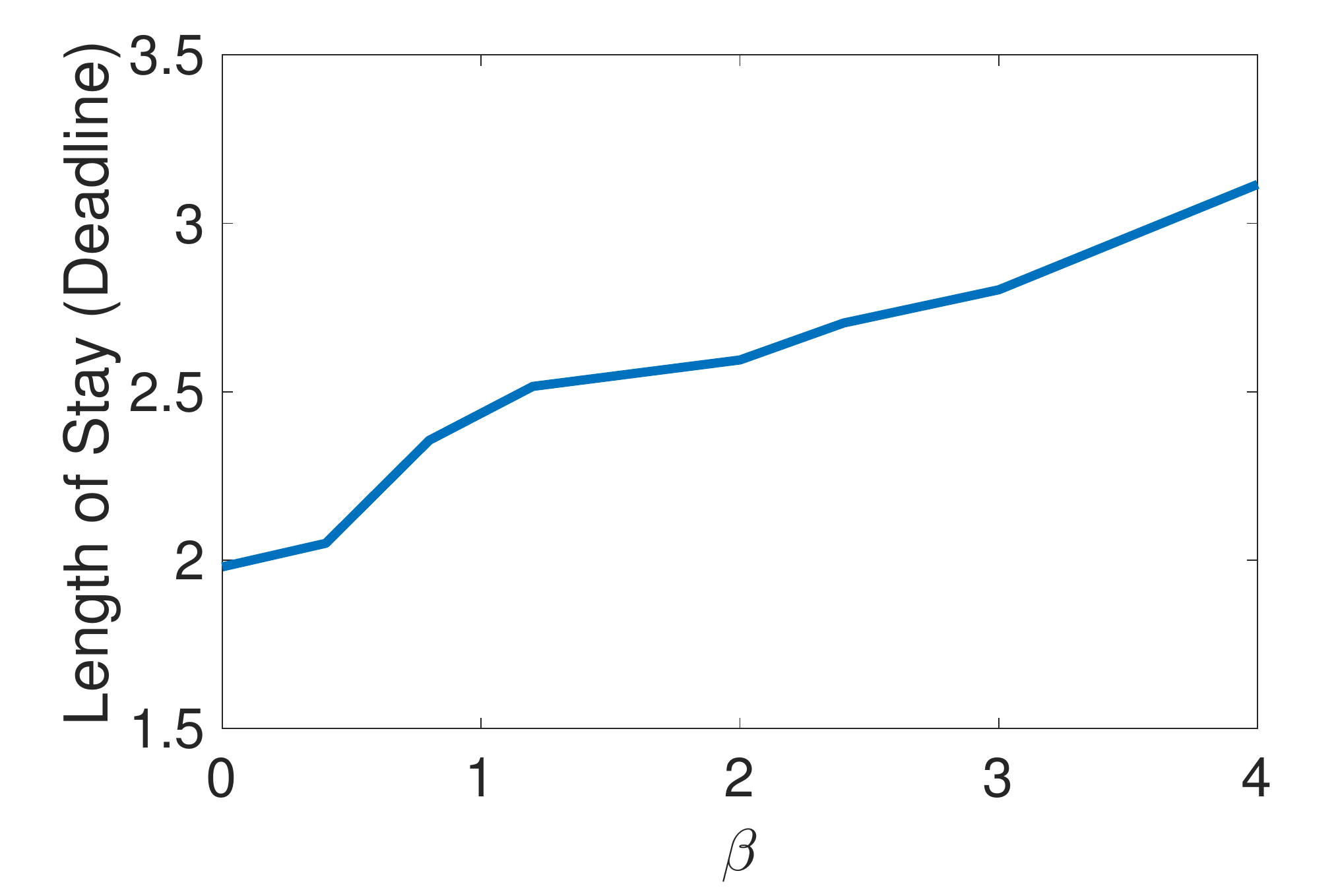}
		\vspace{-0.3in}
		\caption{Variation of the average time spent per EV with $\beta$.}
		\label{fig:deadline}
		\vspace{-0.2in}
	\end{minipage}
\end{figure*}

The total  surpluses of the users  decreases with increase in $\beta$ (Fig.~\ref{fig:consumer_surplus}) as the user pays larger price when $\beta$ increases. The user surplus is maximum at $\beta=0$.  The decrement of total users' surplus is not significant with $\beta$ for $\beta<1.6$. However, as $\beta>1.6$,   it decreases rapidly. For $\beta<1.6$, the number of users served does not decrease much with $\beta$. Hence, the total users' surpluses decrease slowly.

As $\beta$ increases the profit increases initially (Fig.~\ref{fig:profit}). However, as $\beta>3$, the number of users served decreases rapidly, hence, the profit also drops.  

At high values of $\beta$ both users' surpluses and the profit decrease significantly. Low values of $\beta$ give high users' surpluses, however, the profit is low.  $\beta\in [0.8, 1.6]$ is the best candidate for the balance between profit and users' surpluses. 
\subsubsection{Effect on the average deadline}
Our analysis shows that users spend more time in the charging station with the increase in $\beta$  (Fig.~\ref{fig:deadline}). As $\beta$ increases, the users which have preferences for lower deadlines have to pay more; since the cost of fulfilling lower deadline contracts is  high. Hence, those users are reluctant to accept the contract. Thus, the accepted users spend more time in the charging station. Though the increment of the average time spent by an EV is not exponential with $\beta$. The average time spent by an EV is  $2.5$ hours  for $\beta=1.2$ which is in accordance with the average time spent by an EV  \cite{gov}.

\begin{figure*}
	\begin{minipage}{0.19\linewidth}
		\includegraphics[trim=.3in 3in 0.8in 3.4in, clip,width=\textwidth]{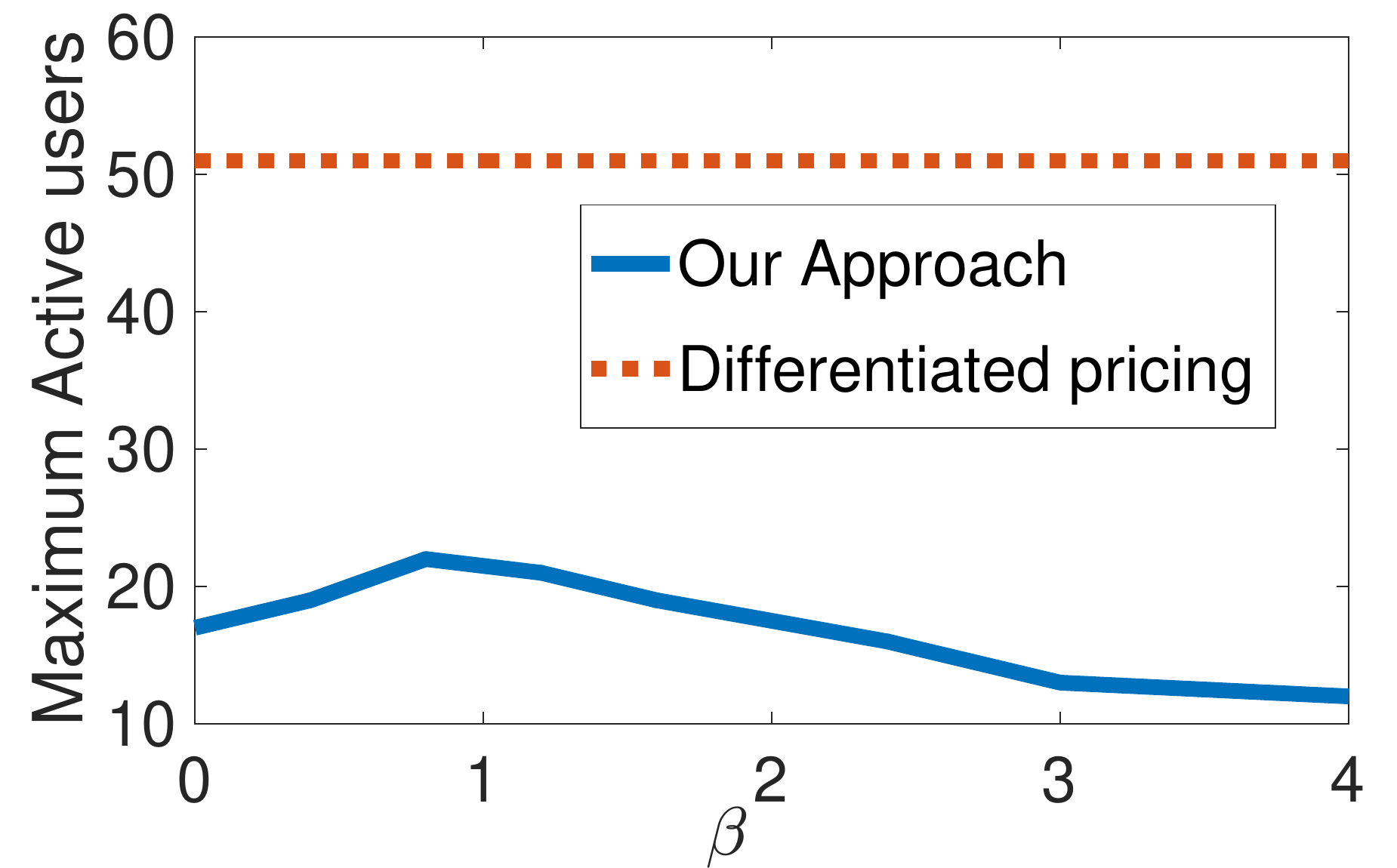}
		\vspace{-0.3in}
		\caption{Variation of the maximum of the number of active users  with $\beta$ and comparison with the differentiated pricing scheme proposed in \cite{bitar,bitar2} (in dotted line).}
		\label{fig:active_max}
		\vspace{-0.4in}
	\end{minipage}\hfill
	\begin{minipage}{0.23\linewidth}
		\includegraphics[trim=0in 0in 0.6in 0in,width=0.99\textwidth]{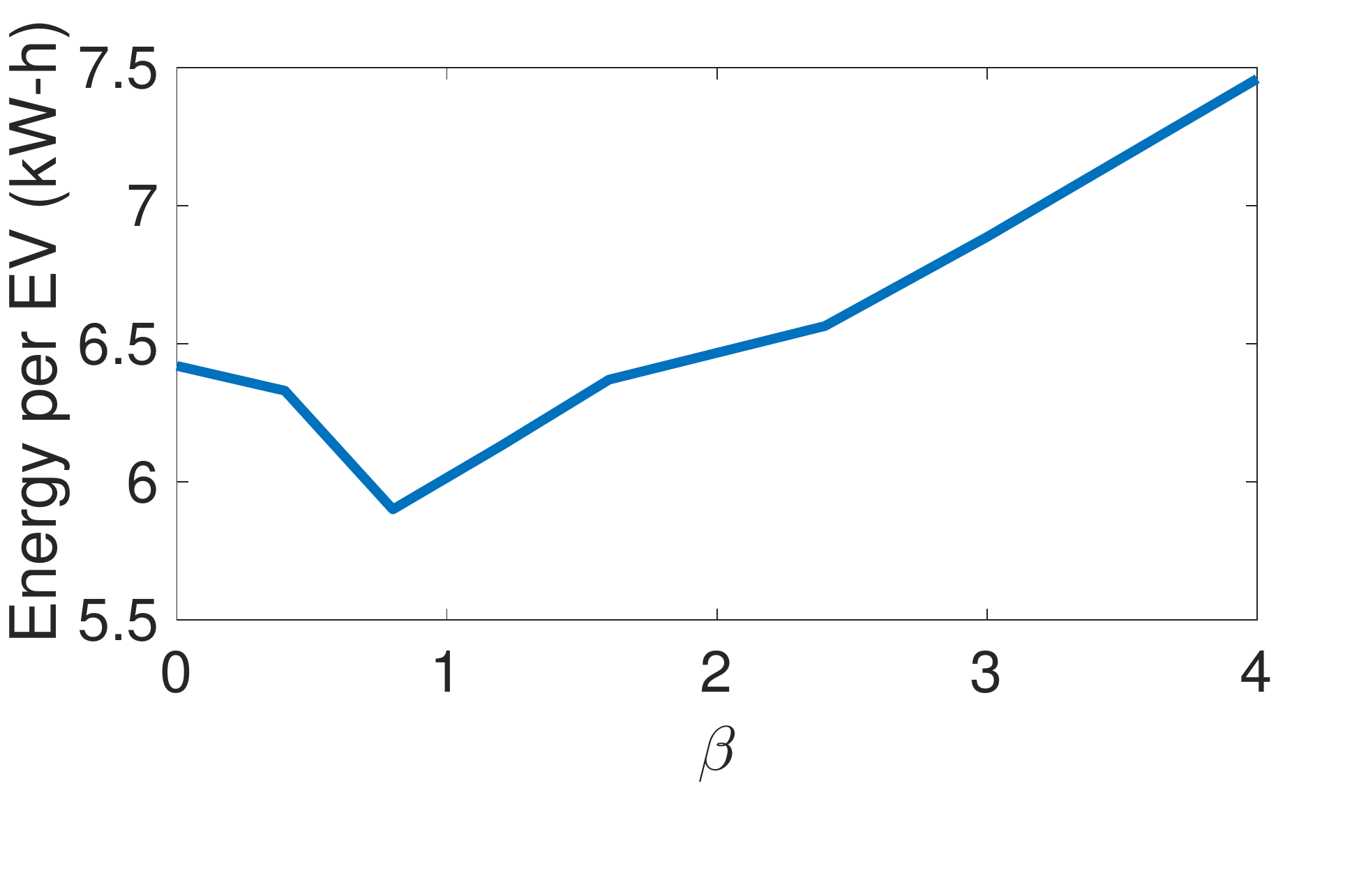}
		\vspace{-0.3in}
		\caption{Variation of the average energy consumed per EV with $\beta$.}
		\label{fig:energy_mean}
		\vspace{-0.4in}
	\end{minipage}
	\begin{minipage}{0.23\linewidth}
		\includegraphics[trim=0in 0in 0.6in 0.2in,width=0.99\textwidth]{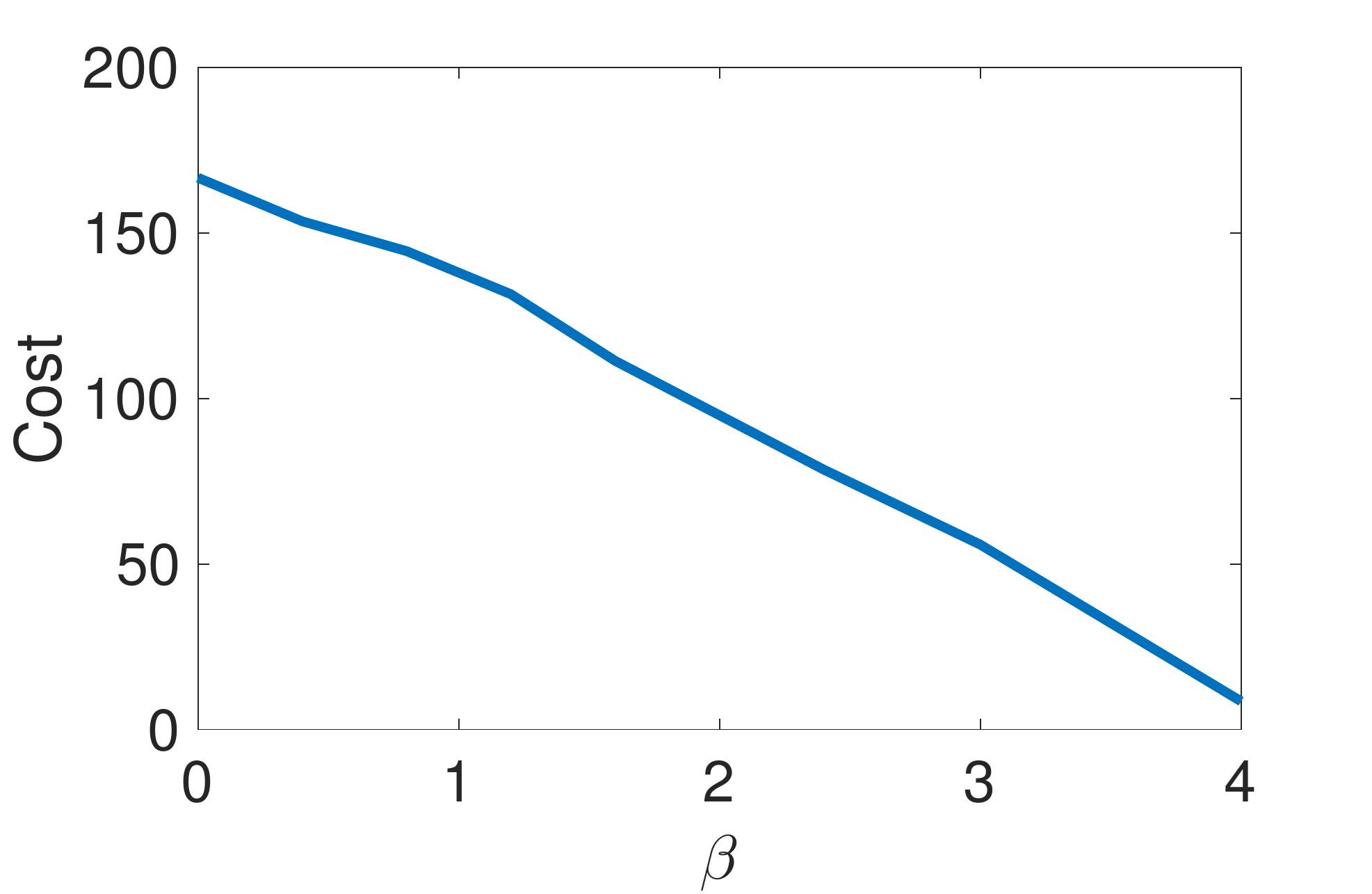}
		\vspace{-0.3in}
		\caption{Variation of the cost incurred by the charging station with $\beta$.}
		\label{fig:cost}
		\vspace{-0.4in}
	\end{minipage}\hfill
	\begin{minipage}{0.3\linewidth}
		\includegraphics[trim=.6in 0in 0.3in 0in,width=0.99\textwidth]{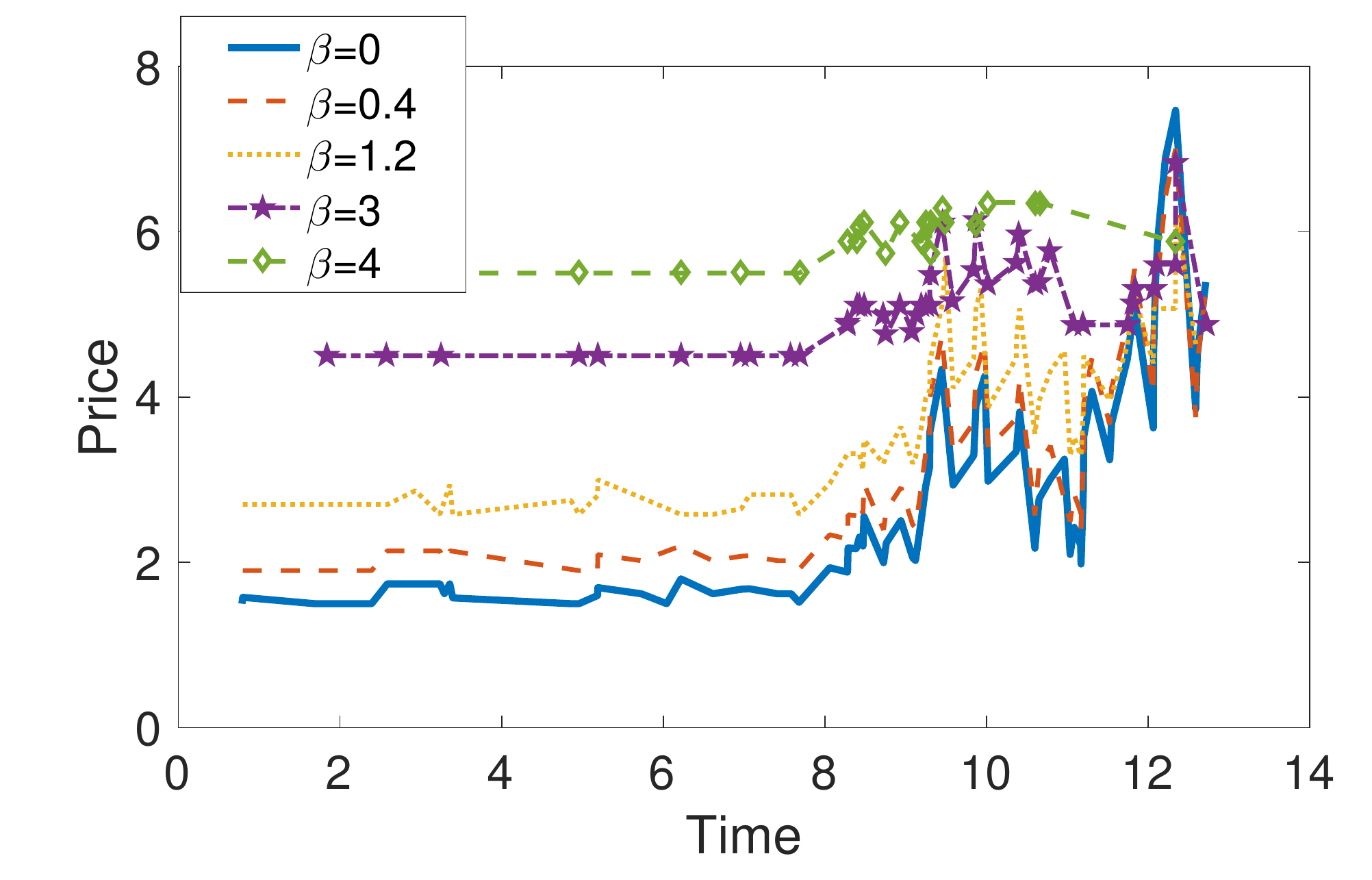}
		\vspace{-0.3in}
		\caption{Variation of the prices set at different times with $\beta$.}
		\label{fig:price}
		\vspace{-0.3in}
	\end{minipage}
\end{figure*}

\subsubsection{Effect on the maximum number of active users}
Since the average time spent by users in the charging station increases with $\beta$ and the number of admitted users are almost the same for $\beta\leq 1.2$, hence the number of active users increases initially as $\beta$ increases (Fig.~\ref{fig:active_max}). Though the maximum never reaches beyond $22$ for any value of $\beta$. However, when $\beta>1.2$, the number of active users decreases with $\beta$.
\subsubsection{Advantages of our proposed mechanism}
Fig.~\ref{fig:active_max} shows that  our pricing algorithm requires less charging spots compared to the differentiated  pricing mechanism \cite{bitar, bitar2} closest to our proposed approach. Similar to \cite{bitar, bitar2}  the users  select the amount of energy to be consumed for each time period based on the price set by the charging station. We assume that the user will not charge beyond the preferred deadline and before the arrival time. In \cite{bitar,bitar2}  the EVs tend to spend more time as it reduces the cost\footnote{An EV is removed when it is fully charged.} and thus, the maximum of the number of EVs present at any time is also higher (Fig.~\ref{fig:active_max}) compared to our proposed mechanism.\footnote{We assume that the EV is removed after its reported deadline. When the deadline is over,  the EV is moved to a parking spot without any charging facility.} In our proposed mechanism, the charging station controls the time spent by an EV through pricing and results into lower charging spots.

\subsubsection{Effect on the average energy}
As $\beta$ increases  users only with higher utilities should accept the contracts. Thus,  the average charging amount for each EV should increase with $\beta$. However, Fig.~\ref{fig:energy_mean} shows that for $\beta\leq 0.8$, the average energy consumed by each EV decreases with the increase in $\beta$. The apparent anomaly is due to the fact that  the users with higher demand but with smaller deadline preferences, may have to pay more because of the increase in the price to fulfill the contract as $\beta$ increases. Hence, such users will not accept the offers which result into initial decrement of the average energy consumption with the increase in $\beta$. However, as $\beta$ becomes large, only the users with higher demand accept the offers, hence, the average energy consumption increases. However, the increment is only linear. In fact for $\beta=2$, the average energy consumption per EV is around $6.9$ kW-h. 

\subsubsection{Effect on the Cost of the EV charging station}
The cost of the EV charging station decreases with the increase in $\beta$ (Fig.~\ref{fig:cost}). Since the time spent by users increases and  thus, the demand of the users can be met through renewable energies. The charging station buys a lower amount of conventional energies which results in lower cost for the charging station.  When $\beta\leq 1.6$, the number of admitted users decreases {\em sub-linearly}, still the cost decreases {\em linearly}. Hence, the FERC will prefer this setting as it decreases the cost without decreasing the admitted users much.  

\subsubsection{Effect on the price selected by the charging station}
The price is higher during the peak period when the arrival rates is higher and the time-of-use price is high (Fig.~\ref{fig:price}). Hence, the pricing mechanism is consistent with the FERC's objective of selecting higher prices during the peak time to flatten the demand curve. A new price is selected when an EV arrives. As $\beta$ decreases the admitted users is higher, hence the price variation is also higher as $\beta$ decreases. Also note that when the  number of active users is large,  serving additional user can be significant and thus, the price is also high.

\begin{figure}
	\includegraphics[trim=.6in 0in 1.1in 0in, clip, width=.95\textwidth]{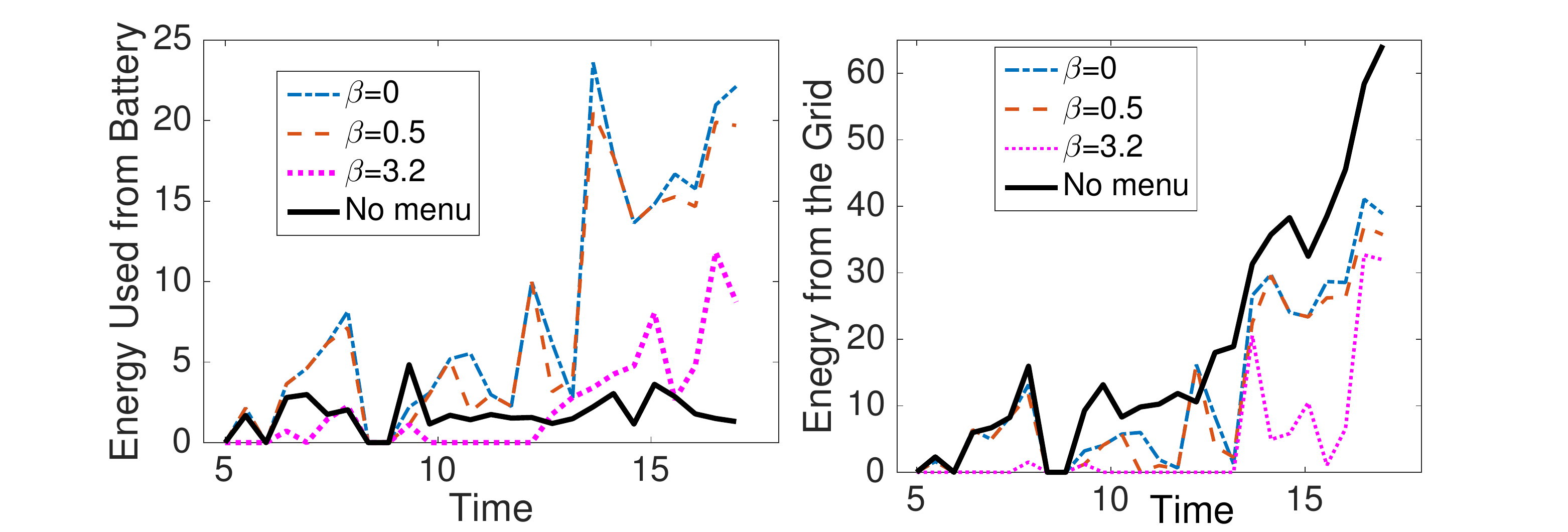}
	\vspace{-0.1in}
	\caption{\small Left-hand figure shows the amount of energy drawn from the battery of the charging station at various times for different values of $\beta$. The right-hand figure shows the amount of energy drawn from the grid at various times for different values of $\beta$.}
	\label{fig:gridvsbeta}
\end{figure}
\subsubsection{Impact on the energy drawn from the grid and the storage of the charging station}
Fig.~\ref{fig:gridvsbeta} shows that as $\beta$ increases the energy bought from the grid decreases. This is because the number of accepted users decreases with $\beta$. The energy used from the battery also decreases as $\beta$ increases. Note that by selecting a $\beta$ the charging station can also limit the peak energy consumption from the grid.  Fig. ~\ref{fig:gridvsbeta} also shows that when no menu based pricing is applied i.e., the EVs are charged as soon as they arrive, then, the peak energy consumption from the grid is very high. Even $\beta=0$ lowers the energy consumption from the grid significantly. The energy used from the battery of the charging station is also low when there is no menu-based pricing. {\em The above shows the usefulness of the menu-based pricing in reducing the peak-energy consumption and efficient use of the renewable energy.}

\begin{figure}
	\includegraphics[trim=.6in 0in 1.1in 0in, clip, width=.95\textwidth]{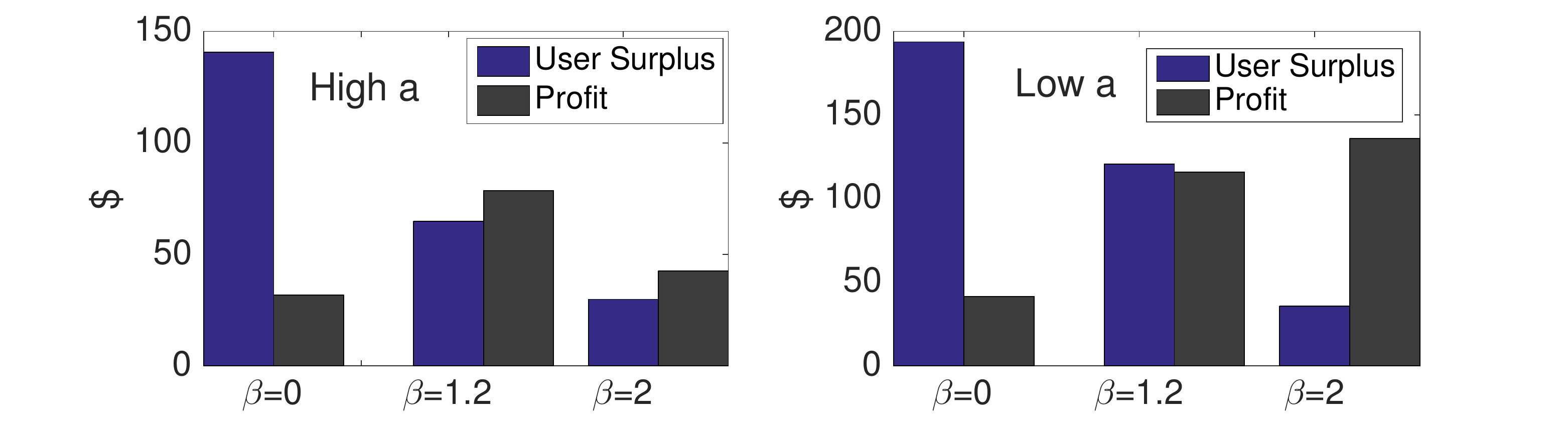}
	\vspace{-0.1in}
	\caption{\small In the left-hand figure, we consider $a\sim\mathcal{U}[1/40,1/16]$. In the right-hand figure, we consider $a\sim \mathcal{U}[1/10,1/4]$. We show the consumer surplus and profit of the charging station for $\beta=0, 1.2,2$.}
	\label{fig:highvslowa}
	\vspace{-0.2in}
\end{figure}
\subsubsection{Impact of $a$}
Fig.~\ref{fig:highvslowa} shows that as $a$ increases, the profit and the user's surpluses both decrease. Note  that as $a$ increases, the utility decreases, and the preferred energy $\dfrac{b}{2a}$ also decreases, hence, the profit and the user's surplus both decrease.
\section{Conclusions and Future Works}
We propose an online menu-based pricing mechanism for EV charging. Specifically, we consider that the charging station will offer a price to each arriving user for a plethora of options; the user selects either one of them or rejects all. We show that there exists a {\em prior-free} pricing strategy which maximizes the ex-post social welfare. We characterize the maximum possible profit that the charging station can get while maximizing the {\em ex-post} social welfare. The charging station only needs to know the lower end-points of the utilities to implement the pricing strategy. The profit  increases if the renewable energy penetration increases or the storage capacity of the charging station increases. However, the increment is bounded. 

The charging station can not simultaneously maximize the profit and the ex-post social welfare unless it is {\em clairvoyant}. We propose a fixed profit pricing scheme which provides a fixed profit to the charging station. The above can also maximize the expected profit of the charging station under some assumptions which frequently arise in practice. Numerical evaluation suggests that the menu-based pricing scheme can reduce the peak-demand and utilize the limited number of charging spots more efficiently compared to the baseline approaches.

Following this work, we have considered the case where the EVs can inject energies by discharging via a Vehicle-to-Grid (V2G) service which can enhance the profits of the charging station \cite{2016arXiv161200106G}.

We considered that the EV charging station is myopic which does not consider the future arrival process while selecting an optimal price for an incoming EV.  In future we consider the case where the charging station  knows the statistics of the future arrival process of the EVs and selects price accordingly. We also considered that the charging station has the only type of charger (either fast or slow), the characterization of prices when the charging station selects prices for different chargers is also left for the future. Considering stochastic pattern of  energy harvesting is an important next step. Finally, the consideration of the  multiple charging stations which set prices in a competitive manner also constitutes a future research direction.

\bibliographystyle{IEEEtran}
\bibliography{ev_charging}

\begin{thebibliography}{10}
\providecommand{\url}[1]{#1}
\csname url@samestyle\endcsname
\providecommand{\newblock}{\relax}
\providecommand{\bibinfo}[2]{#2}
\providecommand{\BIBentrySTDinterwordspacing}{\spaceskip=0pt\relax}
\providecommand{\BIBentryALTinterwordstretchfactor}{4}
\providecommand{\BIBentryALTinterwordspacing}{\spaceskip=\fontdimen2\font plus
\BIBentryALTinterwordstretchfactor\fontdimen3\font minus
  \fontdimen4\font\relax}
\providecommand{\BIBforeignlanguage}[2]{{%
\expandafter\ifx\csname l@#1\endcsname\relax
\typeout{** WARNING: IEEEtran.bst: No hyphenation pattern has been}%
\typeout{** loaded for the language `#1'. Using the pattern for}%
\typeout{** the default language instead.}%
\else
\language=\csname l@#1\endcsname
\fi
#2}}
\providecommand{\BIBdecl}{\relax}
\BIBdecl

\bibitem{icc_veh17}
A.~Ghosh and V.~Aggarwal, ``Control of charging of electric vehicles through
  menu-based pricing under uncertainty,'' in \emph{2017 IEEE International
  Conference on Communications (ICC)}, May 2017.

\bibitem{tech_ev}
------, ``Control of electric vehicles through menu based pricing,''
  \emph{CoRR}, vol. abs/1609.09037, 2016.

\bibitem{soltani}
N.~Y. Soltani, S.~J. Kim, and G.~B. Giannakis, ``Real-time load elasticity
  tracking and pricing for electric vehicle charging,'' \emph{IEEE Transactions
  on Smart Grid}, vol.~6, no.~3, pp. 1303--1313, May 2015.

\bibitem{oren}
Z.~Liu, Q.~Wu, S.~Oren, S.~Huang, R.~Li, and L.~Cheng, ``Distribution
  locational marginal pricing for optimal electric vehicle charging through
  chance constrained mixed-integer programming,'' \emph{IEEE Transactions on
  Smart Grid}, vol.~PP, no.~99, pp. 1--1, 2016.

\bibitem{javidi}
M.~Alizadeh, H.~T. Wai, M.~Chowdhury, A.~Goldsmith, A.~Scaglione, and
  T.~Javidi, ``Optimal pricing to manage electric vehicles in coupled power and
  transportation networks,'' \emph{IEEE Transactions on Control of Network
  Systems}, vol.~PP, no.~99, pp. 1--1, 2016.

\bibitem{tvt}
S.~G. Yoon, Y.~J. Choi, J.~K. Park, and S.~Bahk, ``Stackelberg-game-based
  demand response for at-home electric vehicle charging,'' \emph{IEEE
  Transactions on Vehicular Technology}, vol.~65, no.~6, pp. 4172--4184, June
  2016.

\bibitem{kar}
S.~Bhattacharya, K.~Kar, J.~H. Chow, and A.~Gupta, ``Extended second price
  auctions with elastic supply for pev charging in the smart grid,'' \emph{IEEE
  Transactions on Smart Grid}, vol.~7, no.~4, pp. 2082--2093, July 2016.

\bibitem{zou}
S.~Zou, Z.~Ma, and X.~Liu, ``Distributed efficient charging coordinations for
  electric vehicles under progressive second price auction mechanism,'' in
  \emph{52nd IEEE Conference on Decision and Control}, Dec 2013, pp. 550--555.

\bibitem{sen}
C.~Joe-Wong, S.~Sen, S.~Ha, and M.~Chiang, ``Optimized day-ahead pricing for
  smart grids with device-specific scheduling flexibility,'' \emph{IEEE Journal
  on Selected Areas in Communications}, vol.~30, no.~6, pp. 1075--1085, July
  2012.

\bibitem{low2}
N.~Li, L.~Chen, and S.~H. Low, ``Optimal demand response based on utility
  maximization in power networks,'' in \emph{2011 IEEE Power and Energy Society
  General Meeting}, July 2011, pp. 1--8.

\bibitem{parkes}
\BIBentryALTinterwordspacing
E.~H. Gerding, V.~Robu, S.~Stein, D.~C. Parkes, A.~Rogers, and N.~R. Jennings,
  ``Online mechanism design for electric vehicle charging,'' in \emph{The 10th
  International Conference on Autonomous Agents and Multiagent Systems - Volume
  2}, ser. AAMAS '11.\hskip 1em plus 0.5em minus 0.4em\relax Richland, SC:
  International Foundation for Autonomous Agents and Multiagent Systems, 2011,
  pp. 811--818. [Online]. Available:
  \url{http://dl.acm.org/citation.cfm?id=2031678.2031733}
\BIBentrySTDinterwordspacing

\bibitem{tong}
S.~Chen, Y.~Ji, and L.~Tong, ``Large scale charging of electric vehicles,'' in
  \emph{2012 IEEE Power and Energy Society General Meeting}, July 2012, pp.
  1--9.

\bibitem{qhuang}
Q.~Huang, Q.~S. Jia, Z.~Qiu, X.~Guan, and G.~Deconinck, ``Matching ev charging
  load with uncertain wind: A simulation-based policy improvement approach,''
  \emph{IEEE Transactions on Smart Grid}, vol.~6, no.~3, pp. 1425--1433, May
  2015.

\bibitem{xu_cdc}
Y.~Xu and F.~Pan, ``Scheduling for charging plug-in hybrid electric vehicles,''
  in \emph{2012 IEEE 51st IEEE Conference on Decision and Control (CDC)}, Dec
  2012, pp. 2495--2501.

\bibitem{yu_allerton}
Z.~Yu, Y.~Xu, and L.~Tong, ``Large scale charging of electric vehicles: A
  multi-armed bandit approach,'' in \emph{53rd Annual Allerton Conference on
  Communication, Control, and Computing (Allerton)}, 2015, pp. 389--395.

\bibitem{bitar}
E.~Bitar and S.~Low, ``Deadline differentiated pricing of deferrable electric
  power service,'' in \emph{2012 IEEE 51st IEEE Conference on Decision and
  Control (CDC)}, Dec 2012, pp. 4991--4997.

\bibitem{bitar2}
E.~Bitar and Y.~Xu, ``Deadline differentiated pricing of deferrable electric
  loads,'' \emph{IEEE Transactions on Smart Grid}, vol.~8, no.~1, pp. 13--25,
  Jan 2017.

\bibitem{nayyar}
A.~Nayyar, M.~Negrete-Pincetic, K.~Poolla, and P.~Varaiya,
  ``Duration-differentiated energy services with a continuum of loads,''
  \emph{IEEE Transactions on Control of Network Systems}, vol.~3, no.~2, pp.
  182--191, June 2016.

\bibitem{Salah2016}
F.~Salah and C.~M. Flath, ``Deadline differentiated pricing in practice:
  marketing ev charging in car parks,'' \emph{Computer Science - Research and
  Development}, vol.~31, no.~1, pp. 33--40, 2016.

\bibitem{Hinterhuber}
\BIBentryALTinterwordspacing
A.~Hinterhuber, ``Towards value-based pricing—an integrative framework for
  decision making,'' \emph{Industrial Marketing Management}, vol.~33, no.~8,
  pp. 765 -- 778, 2004. [Online]. Available:
  \url{http://www.sciencedirect.com/science/article/pii/S0019850103001524}
\BIBentrySTDinterwordspacing

\bibitem{personalized_pricing}
M.~H. Yaghmaee, M.~S. Kouhi, and A.~L. Garcia, ``Personalized pricing: A new
  approach for dynamic pricing in the smart grid,'' in \emph{2016 IEEE Smart
  Energy Grid Engineering (SEGE)}, Aug 2016, pp. 46--51.

\bibitem{quadratic}
W.~Shuai, P.~Maill{\'{e}}, and A.~Pelov, ``Charging electric vehicles in the
  smart city: A survey of economy-driven approaches,'' \emph{IEEE Transactions
  on Intelligent Transportation Systems}, vol.~17, no.~8, pp. 2089--2106, Aug
  2016.

\bibitem{low}
L.~Gan, U.~Topcu, and S.~H. Low, ``Optimal decentralized protocol for electric
  vehicle charging,'' \emph{IEEE Transactions on Power Systems}, vol.~28,
  no.~2, pp. 940--951, May 2013.

\bibitem{gov}
``Evaluating electric vehicle charging impacts and customer charging behaviors-
  experiences from six smart grid investment grant projects,''
  \url{https://www.smartgrid.gov/files/B3_revised_master-12-17-2014_report.pdf},
  December, 2014.

\bibitem{2016arXiv161200106G}
A.~{Ghosh} and V.~{Aggarwal}, ``{Menu-Based Pricing for Charging of Electric
  Vehicles with Vehicle-to-Grid Service},'' \emph{ArXiv e-prints}, vol.
  1612.00106, Nov. 2016.

\end{thebibliography}

\begin{IEEEbiography}[{\includegraphics[width=1in,height=1.25in,clip,keepaspectratio]{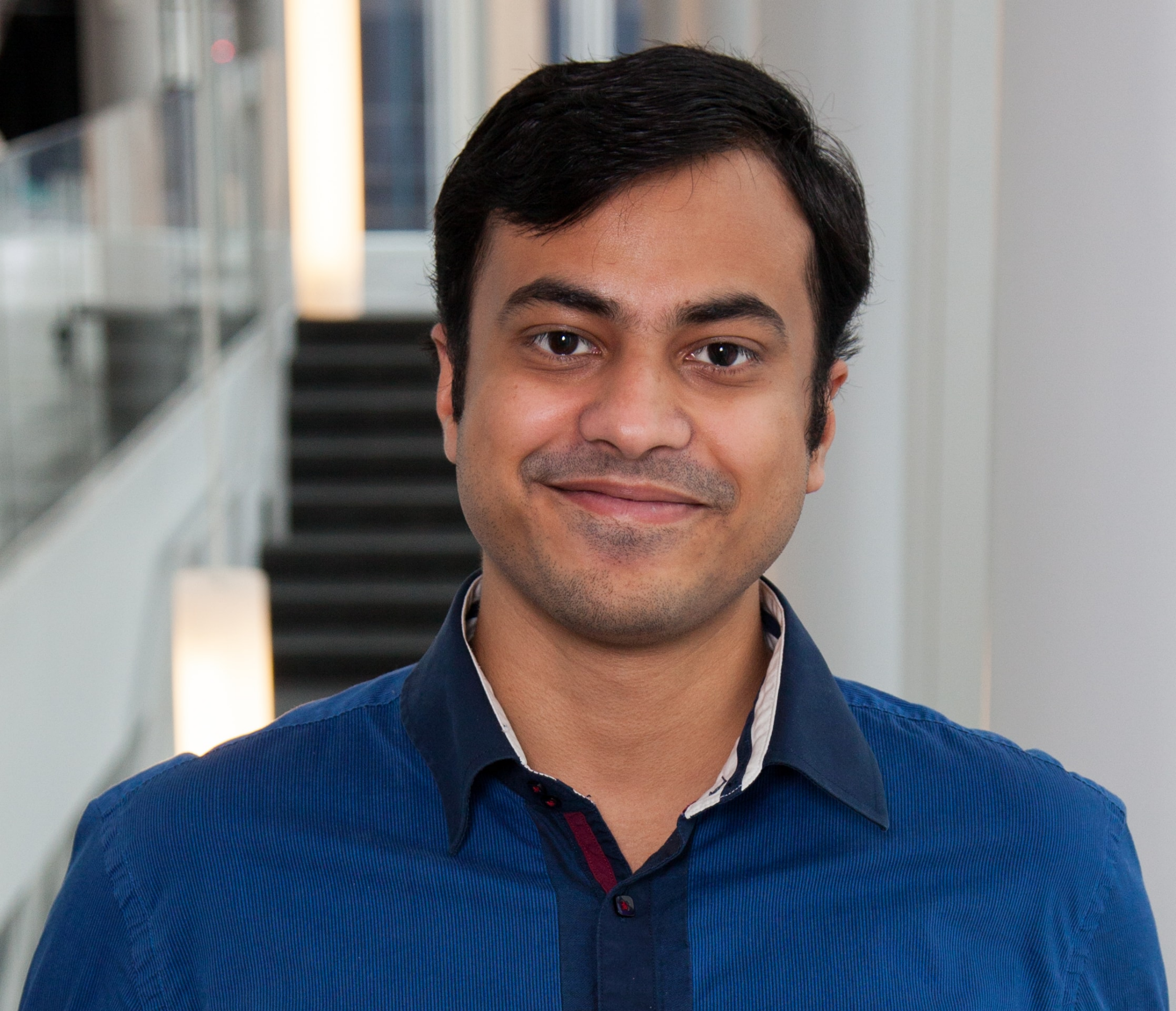}}]{Arnob Ghosh} received the B.E. degree in Electronics and Telecommunications from Jadavpur University, Kolkata, India in 2011; the PhD and M.S. degree in Electrical Engineering from the University of Pennsylvania in 2016 and 2013 respectively. He is currently a Postdoctoral Research Associate in the Industrial Engg. Department of Purdue University. His research interests are smart grid, economic aspects of spectrum sharing in wireless network, and game theoretic approach for resource allocations in wireless network. He served as a reviewer in {\it IEEE Transactions on Signal Processing}, {\it ACM/IEEE Transactions on Networking}, and {\it IEEE Transactions on Communication}. 
\end{IEEEbiography}
\begin{IEEEbiography}[{\includegraphics[width=1in,height=1.25in,clip,keepaspectratio]{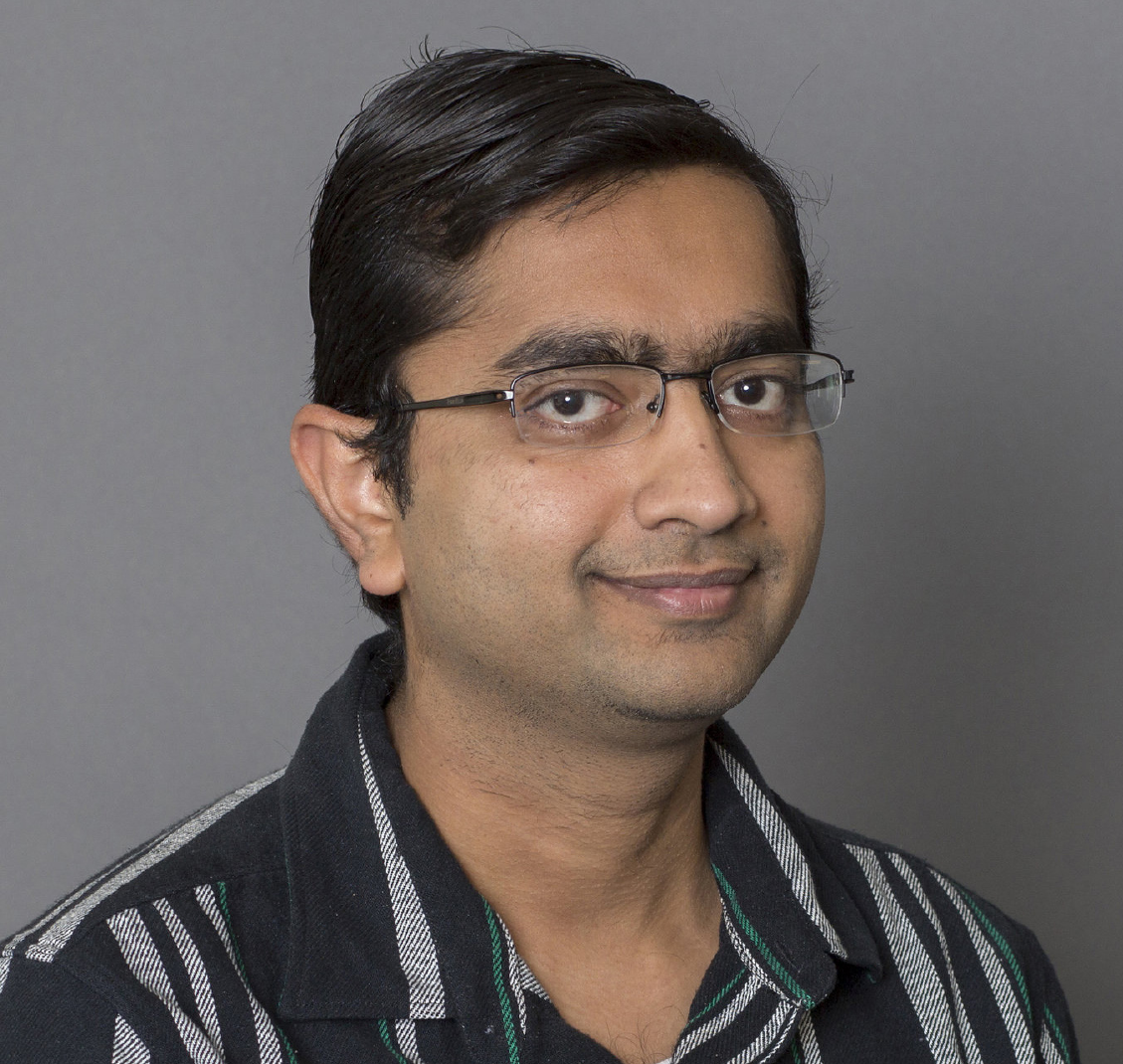}}]{Vaneet Aggarwal (S'08 - M'11 - SM’15)}
	received the B.Tech. degree in 2005 from the Indian Institute of Technology, Kanpur, India, and the M.A. and Ph.D. degrees in 2007 and 2010, respectively from Princeton University, Princeton, NJ, USA, all in Electrical Engineering.
	
	He is currently an Assistant Professor at Purdue University, West Lafayette, IN. Prior to this, he was a Senior Member of Technical Staff Research at AT\&T Labs-Research, NJ, and an Adjunct Assistant Professor at Columbia University, NY. His current research interests are in communications and networking, machine learning, and smart grids. Dr. Aggarwal was the recipient of Princeton University's Porter Ogden Jacobus Honorific Fellowship in 2009. In addition, he received AT\&T Key Contributor award in 2013, AT\&T Vice President Excellence Award in 2012, and AT\&T Senior Vice President Excellence Award in 2014. He is serving on the editorial board of the {\it IEEE Transactions on Communications} and the {\it IEEE Transactions on Green Communications and Networking.}
\end{IEEEbiography}

\appendix
\subsection{Proof of Theorem~\ref{thm:pricestrategy}}
Note that the user $k$  only selects the menu which fetches the highest payoff. Hence, the payoff of user $k$ is
\begin{align}
(\max_{l,t}\{ u_{k,l,t}-p_{k,l,t}\})^{+}
\end{align}
Since $p_{k,l,t}=v_{l,t}-v_{-k}$, thus, the usersurplus is
\begin{align}
(\max_{l,t}(u_{k,l,t}-v_{l,t}+v_{-k}))^{+}
\end{align}
The charging station's profit is $p_{k,l,t}-v_{l,t}$ if the user selects the menu $(l,t)$ and $-v_{-k}$ if the user does not select any price option. Since $p_{k,l,t}=v_{l,t}-v_{-k}$, thus, the charging station's profit is $-v_{-k}$ irrespective of the decision of the user. Hence, the social welfare is
\begin{align}
\max\{ \max_{l,t}(u_{k,l,t}-v_{l,t}),v_{-k}\}
\end{align}

Hence, under the pricing strategy, the value of the social welfare is the same as in Theorem~\ref{thm:max_socialwelfare}.\qed

\subsection{Proof of Theorem~\ref{thm:profit_max}}
First, we show that the pricing strategy maximizes the profit of the charging station. 

From (\ref{eq:profit_maximum}), the highest possible profit of the charging station is
\begin{align}
\max\{u_{k,l^{*},t^{*}}-v_{l^{*},t^{*}}, -v_{-k}\}
\end{align}
If the price is selected as stated in the theorem, then, if $u_{k,l^{*},t^{*}}-v_{l^{*},t^{*}}<-v_{-k}$, then, the profit of the charging station is $-v_{-k}$.

On the other hand, if $u_{k,l^{*},t^{*}}-v_{l^{*},t^{*}}\geq -v_{-k}$, then the price $p_{k,l,t}$ to user $k$ is
\begin{align}
p_{k,l,t}=v_{l,t}+u_{k,l^{*},t^{*}}-v_{l^{*},t^{*}}
\end{align}
User will select a price $p_{k,l,t}$ if $u_{k,l,t}-p_{k,l,t}\geq u_{k,i,j}-p_{k,i,j}$ for all $i,j$ and $u_{k,l,t}\geq p_{k,l,t}$. 
Since $u_{k,l^{*},t^{*}}-v_{l^{*},t^{*}}\geq u_{k,l,t}-v_{l,t}$, thus, the user will only select the price menu $p_{k,l^{*},t^{*}}$. Note that at only $p_{k,l^{*},t^{*}}$ the user surplus is $0$, at other prices it is less than or equal to $0$.  Thus, the profit of the charging station is $u_{k,l^{*},t^{*}}-v_{l^{*},t^{*}}$.

Now, we show that such a pricing scheme also maximizes the social welfare. 

Note that in the above pricing strategy, the user's surplus is always $0$. On the other hand, the profit that the charging station makes is $\max\{\max_{l,t}(u_{k,l,t}-v_{l,t}), -v_{-k}\}$. Hence, the pricing strategy obtains the optimal value of social welfare by Theorem~\ref{thm:max_socialwelfare}.\qed

\subsection{Proof of Theorem~\ref{thm:profitmax_uncertainty}}
If $\max_{i,j}\{L_{k,i,j}-v_{i,j}+v_{-k}\}\leq 0$, then the pricing strategy is the same as in Theorem~\ref{thm:price_expected} which we already proved to be ex-post social welfare maximizer. Hence, we consider the case when $\max_{i,j}\{L_{k,i,j}-v_{i,j}+v_{-k}\}>0$.

Now, we show that for every possible realization  $u_{k,l,t}$ such a pricing strategy will maximize the social welfare. Note that $\max_{l,t}\{u_{k,l,t}-v_{l,t}+v_{-k}\}>0$ since $\max_{l,t}\{L_{k,l,t}-v_{l,t}+v_{-k}\}>0$ and $L_{k,l,t}$ is the lowest end-point of the distribution of $U_{k,l,t}$. Hence by Theorem~\ref{thm:max_socialwelfare} the maximum value of the ex-post social welfare is $\max_{l,t}\{u_{k,l,t}-v_{l,t}\}$. Now, we show that the pricing strategy defined in (\ref{eq:pricesocialmax}) will give rise the above optimal value of social welfare.

Note that  the user's surplus is 
\begin{align}
(\max_{i,j}\{u_{k,i,j}-p_{k,i,j}\})^{+}
\end{align}

First, we show that $\max_{i,j}\{u_{k,i,j}-p_{k,i,j}\}\geq 0$ if $\max_{i,j}\{L_{k,i,j}-v_{i,j}+v_{-k}\}>0$ and $p_{k,i,j}$ is given by (\ref{eq:pricesocialmax}). Suppose not, i.e. $\max_{i,j}\{u_{k,i,j}-p_{k,i,j}\}< 0$. Let $(l^{*},t^{*})=\text{argmax}_{i,j}\{L_{k,i,j}-v_{i,j}+v_{-k}\}$. Then, $p_{k,l^{*},t^{*}}=L_{k,l^{*},t^{*}}$. Since $u_{k,l^{*},t^{*}}\geq L_{k,l^{*},t^{*}}$, hence, $u_{k,l^{*},t^{*}}-p_{k,l^{*},t^{*}}>0$ which leads to a contradiction. Thus, $\max_{i,j}\{u_{k,i,j}-p_{k,i,j}\}>0$.

Thus, the user surplus is $\max_{i,j}\{u_{k,i,j}-p_{k,i,j}\}$. Since $p_{k,l,t}=v_{l,t}-v_{-k}+r$ where $r=\max_{i,j}\{L_{k,i,j}-v_{i,j}+v_{-k}\}$ is constant and independent of index $l$ and $t$. Hence, the user will select the price menu $p_{k,l^{*},t^{*}}$ such that $(l^{*},t^{*})=\text{argmax}_{l,t}\{u_{k,l,t}-v_{l,t}\}$ . The profit of the charging station is $p_{k,l^{*},t^{*}}-v_{l^{*},t^{*}}$. Thus, the ex-post social welfare is
\begin{align}
& u_{k,l^{*},t^{*}}-p_{k,l^{*},t^{*}}+p_{k,l^{*},t^{*}}-v_{l^{*},t^{*}}\nonumber\\
& =u_{k,l^{*},t^{*}}-v_{l^{*},t^{*}}\nonumber\\
& =\max_{l,t}\{u_{k,l,t}-v_{l,t}\}
\end{align}
Thus, for each realized values of the utilities $u_{k,l,t}$, the pricing strategy (\ref{eq:pricesocialmax}) provides the maximum social welfare. Hence, the result follows.\qed
\subsection{Proof of Theorem~\ref{thm:approx}}
Suppose that $u_{k,l,t}$ be the realized values of $U_{k,l,t}$. Note that if $\max_{l,t}\{u_{k,l,t}-p_{k,l,t}\}\geq 0$ where $p_{k,l,t}$ is given by (\ref{eq:approx_price}) then the user selects the price menu $p_{k,l^{*},t^{*}}$ where 
\begin{align}
(l^{*},t^{*})=\text{argmax }_{l,t}\{u_{k,l,t}-v_{l,t}\}.
\end{align}
Hence, the social welfare is $\max_{l,t}\{u_{k,l,t}-v_{l,t}\}$ which is the same as the social welfare by Theorem~\ref{thm:max_socialwelfare}.

On the other hand if $\max_{l,t}\{u_{k,l,t}-p_{k,l,t}\}<0 0$, then the social welfare is $-v_{-k}$. 

Hence, by Theorem~\ref{thm:max_socialwelfare}, the only cases where the social welfare is not maximized when $\max_{l,t}\{u_{k,l,t}-p_{k,l,t}\}<0$, however, $\max_{l,t}\{u_{k,l,t}-v_{l,t}+v_{-k}\}\geq 0$. However, by the definition of $\delta(\epsilon)$ such a scenario can only occur with probability $\epsilon$. Hence, the result follows. \qed
\subsection{Proof of Theorem~\ref{thm:aclassutility}}
Suppose the statement is false. Without loss of generality, assume that $p_{k,l,t}=v_{l,t}-v_{-k}+\alpha_{l,t}$ where $\alpha_{l,t}\neq \alpha$ for some $l$ and $t$ achieves a strictly higher expected payoff than the pricing strategy $p_{k,l,t}=v_{l,t}-v_{-k}+\alpha$. 

The expected profit of the charging station for pricing strategy $p_{k,l,t}=v_{l,t}-v_{-k}+\alpha_{l,t}$ is given by
\begin{align}
& \sum_{l=1}^{L}\sum_{t=t_k+1}^{T}(p_{k,l,t}-v_{l,t}+v_{-k})\Pr(R_{l,t})-v_{-k}=\sum_{l=1}^{L}\sum_{t=t_k+1}^{T}\alpha_{l,t}\Pr(R_{l,t})-v_{-k}
\end{align}
Now, we evaluate the expression $\Pr(R_{l,t})$. The user $k$ will select the menu $p_{k,l,t}$ with a positive probability  if $Y_{k,l,t}+X_{k}\geq v_{l,t}-v_{-k}+\alpha_{l,t}$ and for every $(i,j)\neq (l,t)$,
\begin{align}
\alpha_{l,t}-\alpha_{i,j}\leq Y_{k,l,t}-v_{l,t}-Y_{k,i,j}+v_{i,j}
\end{align}
Since $Y_{k,l,t}, Y_{k,i,j}, v_{l,t},v_{i,j}$ are fixed, hence,  the above inequality is either satisfied or not satisfied with probability $1$. More specifically, the user selects the menu $(l,t)$ if $Y_{k,l,t}+X_{k}\geq v_{l,t}-v_{-k}+\alpha_{l,t}$ and
\begin{align}
Y_{k,l,t}-v_{l,t}-\alpha_{l,t}\geq \max_{i,j}(Y_{k,i,j}-v_{i,j}-\alpha_{i,j})
\end{align}
Without loss of generality, assume that $\alpha_{l_1,t_1}$ be the maximum value for which the above inequality is satisfied i.e.
\begin{align}
\alpha_{l_1,t_1}=\max\{\alpha_{l,t}: Y_{k,l,t}-v_{l,t}-\alpha_{l,t}\geq \max_{i,j}(Y_{k,i,j}-v_{i,j}-\alpha_{i,j})\}
\end{align}
The random variable $X_k$ only affects the probability whether $Y_{k,l,t}+X_k\geq v_{l,t}-v_{-k}+\alpha_{l_1,t_1}$ or not. Hence, the charging station's expected  profit is upper bounded by
\begin{align}
\alpha_{l_1,t_1}\Pr(X_k\geq v_{l_1,t_1}-v_{-k}+\alpha_{l_1,t_1}-Y_{k,l_1,t_1})
\end{align}
Note that by the definition of $\alpha$ (Definition~\ref{defn:alpha}),
\begin{align}\label{eq:compare}
& \alpha\max_{l,t}\Pr(Y_{k,l,t}+X_{k}\geq v_{l,t}-v_{-k}+\alpha)\nonumber\\
& \geq \alpha_{l_1,t_1}\Pr(X_k\geq v_{l_1,t_1}-v_{-k}+\alpha_{l_1,t_1}-Y_{k,l_1,t_1})
\end{align}
However by Theorem~\ref{thm:expected_payoff} the expected payoff of the charging station when it selects the price $v_{l,t}-v_{-k}+\alpha$ is given by the expression in the left hand of the expression in (\ref{eq:compare}). Hence, this leads to a contradiction. Thus, the result follows.\qed
\end{document}